\documentclass[12pt]{article}
\usepackage{full page, amsfonts, tikz,
amssymb, amscd, amsmath, graphicx, 
wrapfig, stmaryrd, hyperref}
\usetikzlibrary{shapes}

 \title{{\bf Convergence in 
 conformal field theory}}
 \author{Yi-Zhi Huang}
    \date{\it \small In memory of Professor Gu Chaohao}
    \begin{document}
    \bibliographystyle{alpha}
    \maketitle
\newtheorem{thm}{Theorem}[section]
\newtheorem{defn}[thm]{Definition}
\newtheorem{prop}[thm]{Proposition}
\newtheorem{cor}[thm]{Corollary}
\newtheorem{lemma}[thm]{Lemma}
\newtheorem{rema}[thm]{Remark}
\newtheorem{app}[thm]{Application}
\newtheorem{prob}[thm]{Problem}
\newtheorem{conv}[thm]{Convention}
\newtheorem{conj}[thm]{Conjecture}
\newtheorem{cond}[thm]{Condition}
    \newtheorem{exam}[thm]{Example}
\newtheorem{assum}[thm]{Assumption}
     \newtheorem{nota}[thm]{Notation}
\newcommand{\halmos}{\rule{1ex}{1.4ex}}
\newcommand{\pfbox}{\hspace*{\fill}\mbox{$\halmos$}}

\newcommand{\bkdiamd}{\;\raisebox{.5pt}{\tikz{\node[draw,scale=0.35,diamond,fill=black](){};}}\;}

\newcommand{\nn}{\nonumber \\}

 \newcommand{\res}{\mbox{\rm Res}}
 \newcommand{\ord}{\mbox{\rm ord}}
\renewcommand{\hom}{\mbox{\rm Hom}}
\newcommand{\edo}{\mbox{\rm End}\ }
 \newcommand{\pf}{{\it Proof.}\hspace{2ex}}
 \newcommand{\epf}{\hspace*{\fill}\mbox{$\halmos$}}
 \newcommand{\epfv}{\hspace*{\fill}\mbox{$\halmos$}\vspace{1em}}
 \newcommand{\epfe}{\hspace{2em}\halmos}
\newcommand{\nord}{\mbox{\scriptsize ${\circ\atop\circ}$}}
\newcommand{\wt}{\mbox{\rm wt}\ }
\newcommand{\swt}{\mbox{\rm {\scriptsize wt}}\ }
\newcommand{\lwt}{\mbox{\rm wt}^{L}\;}
\newcommand{\rwt}{\mbox{\rm wt}^{R}\;}
\newcommand{\slwt}{\mbox{\rm {\scriptsize wt}}^{L}\,}
\newcommand{\srwt}{\mbox{\rm {\scriptsize wt}}^{R}\,}
\newcommand{\clr}{\mbox{\rm clr}\ }
\newcommand{\tr}{\mbox{\rm Tr}}
\renewcommand{\H}{\mathbb{H}}
\newcommand{\C}{\mathbb{C}}
\newcommand{\Z}{\mathbb{Z}}
\newcommand{\R}{\mathbb{R}}
\newcommand{\Q}{\mathbb{Q}}
\newcommand{\N}{\mathbb{N}}
\newcommand{\CN}{\mathcal{N}}
\newcommand{\F}{\mathcal{F}}
\newcommand{\I}{\mathcal{I}}
\newcommand{\V}{\mathcal{V}}
\newcommand{\one}{\mathbf{1}}
\newcommand{\BY}{\mathbb{Y}}
\newcommand{\ds}{\displaystyle}

        \newcommand{\ba}{\begin{array}}
        \newcommand{\ea}{\end{array}}
        \newcommand{\be}{\begin{equation}}
        \newcommand{\ee}{\end{equation}}
        \newcommand{\bea}{\begin{eqnarray}}
        \newcommand{\eea}{\end{eqnarray}}
         \newcommand{\lbar}{\bigg\vert}
        \newcommand{\p}{\partial}
        \newcommand{\dps}{\displaystyle}
        \newcommand{\bra}{\langle}
        \newcommand{\ket}{\rangle}

        \newcommand{\ob}{{\rm ob}\,}
        \renewcommand{\hom}{{\rm Hom}}

\newcommand{\A}{\mathcal{A}}
\newcommand{\Y}{\mathcal{Y}}
\newcommand{\End}{\mathrm{End}}
 \newcommand{\rad}{\mbox{\rm rad}}
\renewcommand{\l}{\llfloor}
\renewcommand{\r}{\rrfloor}

\newcommand{\dlt}[3]{#1 ^{-1}\delta \bigg( \frac{#2 #3 }{#1 }\bigg) }

\newcommand{\dlti}[3]{#1 \delta \bigg( \frac{#2 #3 }{#1 ^{-1}}\bigg) }

 \makeatletter
\newlength{\@pxlwd} \newlength{\@rulewd} \newlength{\@pxlht}
\catcode`.=\active \catcode`B=\active \catcode`:=\active
\catcode`|=\active
\def\sprite#1(#2,#3)[#4,#5]{
   \edef\@sprbox{\expandafter\@cdr\string#1\@nil @box}
   \expandafter\newsavebox\csname\@sprbox\endcsname
   \edef#1{\expandafter\usebox\csname\@sprbox\endcsname}
   \expandafter\setbox\csname\@sprbox\endcsname =\hbox\bgroup
   \vbox\bgroup
  \catcode`.=\active\catcode`B=\active\catcode`:=\active\catcode`|=\active
      \@pxlwd=#4 \divide\@pxlwd by #3 \@rulewd=\@pxlwd
      \@pxlht=#5 \divide\@pxlht by #2
      \def .{\hskip \@pxlwd \ignorespaces}
      \def B{\@ifnextchar B{\advance\@rulewd by \@pxlwd}{\vrule
         height \@pxlht width \@rulewd depth 0 pt \@rulewd=\@pxlwd}}
      \def :{\hbox\bgroup\vrule height \@pxlht width 0pt depth
0pt\ignorespaces}
      \def |{\vrule height \@pxlht width 0pt depth 0pt\egroup
         \prevdepth= -1000 pt}
   }
\def\endsprite{\egroup\egroup}
\catcode`.=12 \catcode`B=11 \catcode`:=12 \catcode`|=12\relax
\makeatother

\def\hboxtr{\FormOfHboxtr} 
\sprite{\FormOfHboxtr}(25,25)[0.5 em, 1.2 ex] 

:BBBBBBBBBBBBBBBBBBBBBBBBB | :BB......................B |
:B.B.....................B | :B..B....................B |
:B...B...................B | :B....B..................B |
:B.....B.................B | :B......B................B |
:B.......B...............B | :B........B..............B |
:B.........B.............B | :B..........B............B |
:B...........B...........B | :B............B..........B |
:B.............B.........B | :B..............B........B |
:B...............B.......B | :B................B......B |
:B.................B.....B | :B..................B....B |
:B...................B...B | :B....................B..B |
:B.....................B.B | :B......................BB |
:BBBBBBBBBBBBBBBBBBBBBBBBB |

\endsprite
\def\shboxtr{\FormOfShboxtr} 
\sprite{\FormOfShboxtr}(25,25)[0.3 em, 0.72 ex] 

:BBBBBBBBBBBBBBBBBBBBBBBBB | :BB......................B |
:B.B.....................B | :B..B....................B |
:B...B...................B | :B....B..................B |
:B.....B.................B | :B......B................B |
:B.......B...............B | :B........B..............B |
:B.........B.............B | :B..........B............B |
:B...........B...........B | :B............B..........B |
:B.............B.........B | :B..............B........B |
:B...............B.......B | :B................B......B |
:B.................B.....B | :B..................B....B |
:B...................B...B | :B....................B..B |
:B.....................B.B | :B......................BB |
:BBBBBBBBBBBBBBBBBBBBBBBBB |

\endsprite

\makeatletter
\DeclareFontFamily{U}{tipa}{}
\DeclareFontShape{U}{tipa}{m}{n}{<->tipa10}{}
\newcommand{\arc@char}{{\usefont{U}{tipa}{m}{n}\symbol{62}}}%

\newcommand{\arc}[1]{\mathpalette\arc@arc{#1}}

\newcommand{\arc@arc}[2]{%
  \sbox0{$\m@th#1#2$}%
  \vbox{
    \hbox{\resizebox{\wd0}{\height}{\arc@char}}
    \nointerlineskip
    \box0
  }%
}
\makeatother

\date{}
\maketitle

\begin{abstract}
Convergence and analytic extension
are of fundamental importance in the mathematical construction and study of 
conformal field theory. 
We review some main convergence results, conjectures and problems 
in the construction and
study of conformal field theories using the representation theory of 
vertex operator algebras. We also review the related
analytic extension results, conjectures and problems. 
We discuss the convergence and analytic extensions of products of intertwining operators 
(chiral conformal fields) and of  
$q$-traces and pseudo-$q$-traces of 
products of intertwining operators. We also discuss the convergence results 
related to the sewing operation and the determinant line bundle and 
a higher-genus convergence result.
We then explain conjectures and problems 
on the convergence and analytic extensions
in orbifold conformal field theory and in the cohomology theory 
of vertex operator algebras. 
\end{abstract}



\renewcommand{\theequation}{\thesection.\arabic{equation}}
\renewcommand{\thethm}{\thesection.\arabic{thm}}
\setcounter{equation}{0} \setcounter{thm}{0} 

\section{Introduction}

\renewcommand{\theequation}{\thesection.\arabic{equation}}
\renewcommand{\thethm}{\thesection.\arabic{thm}}
\setcounter{equation}{0} \setcounter{thm}{0}

Quantum filed theory has become an important area in mathematics.
Many mathematical problems are deeply connected to structures studied in 
quantum field theory and have been solved or are expected to be solved
using the ideas and tools developed in the mathematical 
study of quantum field theory.

The most successful quantum field theories are topological ones. Since the
state spaces of these quantum field theories are typically finite dimensional, 
topological field theories usually do not involve convergence
problems. On the other hand, since the state spaces of nontopological quantum 
field theories must be infinite dimensional, convergence problems are often
the most basic ones that we have to solve first. 
For topological quantum field theories that are constructed using some
underlying nontopological quantum field theories, one might also need to 
solve some convergence problems. Such convergence problems are
in fact problems for the underlying nontopological quantum field theories. 

Convergence results and problems in mathematics 
are usually related to existence results
and problems. Such existence results and problems are 
always of fundamental importance in mathematics. For example, 
the first fundamental problem for a differential equation is the existence 
of solutions under suitable conditions. Though one can still 
derive many important results by assuming the existence of a mathematical structure, 
such results in mathematics 
are still conjectures unless the existence is proved. In a problem related to 
analysis, the existence is often a problem about convergence. For example,
the existence of a solution of a differential equation is mostly 
a problem about the convergence of a sequences of 
functions obtained from a suitable iteration procedure. 
In some cases, even though the further development of a problem
was studied using algebraic methods, 
without the existence established by proving suitable convergence results,
the subsequent important results would still be conjectures. 

One main nontopological conformal field theory for which we
have precise definitions and substantial mathematical results
is two-dimensional conformal field theory. 
In this paper, for simplicity, we shall omit 
the words ``two-dimensional'' so that
by conformal field theory, we always mean two-dimensional conformal field theory.
Conformal field theory was studied in physics 
using the approach of operator product expansion starting
from the work of Belavin-Polyakov-Zamolodchikov \cite{BPZ}. 
The fundamental early works of 
Friedan-Shenker \cite{FS}, Verlinde \cite{V}, Segal, Moore-Seiberg \cite{MS}
and others led to some major conjectures on rational conformal field theories.
Around the same time, the representation theory of vertex operator algebras
was developed starting from the works of Borcherds \cite{B} 
and Frenkel-Lepowsky-Meurman \cite{FLM}. The representation theory of 
vertex operator algebras is now one of the main 
mathematical approaches for the construction
and study of conformal field theory. 

In the construction and study of conformal field theories using the 
representation theory of vertex operator algebras, 
almost in every major step,
we have to prove a convergence result. Moreover, 
these convergence results are always obtained 
together with some analytic extension results that
are necessary and important 
for proving further results on conformal field theories. 
To further develop conformal field theory
and apply conformal field theory to solve mathematical problems, 
many convergence and analytic extension 
conjectures and problems still need to be proved and solved. 
Without the proofs of these conjectures and 
the solutions to these problems, we would not and will not be able to solve 
many of the mathematical problems related to conformal field theory. 

In this paper, 
we review some main convergence results, conjectures and problems 
in the construction and
study of two-dimensional 
conformal field theories using the representation theory of 
vertex operator algebras. We also discuss the related 
analytic extension results, conjectures and problems. 
These analytic results and problems are of the fundamental importance
in the study of conformal field theory
and the representation theory of vertex operator algebras. 
We believe that in the future study and applications 
of conformal field theory,
the results and problems on convergence and analytic extensions
will play an even more important role.

 For the definitions and basic properties of vertex operator algebras, 
(lower-bounded and grading-restricted) generalized $V$-modules
 and (logarithmic) intertwining operators, see for example, 
 \cite{FHL}, \cite{tensor1} and \cite{HLZ2}.  For simplicity, in this paper,
 we shall call a logarithmic intertwining operator an intertwining operator.

This paper is organized as follows: In Section 2, as a comparison with 
the main convergence results that we shall discuss in later sections, 
we discuss 
the convergence of products of vertex operators for vertex operator algebras 
and modules. We also briefly discuss the generalizations this type of convergence
in this section. In Sections 3 and 4, we discuss the convergence and analytic 
extensions of products of intertwining operators and
the convergence and analytic extensions of  
$q$-traces and pseudo-$q$-traces of 
products of intertwining operators, respectively. 
In Section 5, we discuss the convergence related to the sewing operation 
of spheres with punctures and local coordinates vanishing at the 
punctures and to the determinant 
line bundle. In Section 6, a higher-genus convergence result 
proved by Gui is discussed. The convergence and 
analytic extension conjectures and problems in 
orbifold conformal field theory and 
the convergence problems in the cohomology theory 
of vertex operator algebras are discussed in Section 6 and 7, respectively.

\paragraph{Acknowledgment} The author gave an online minicourse 
``Convergence in conformal field theory'' based on the material 
in this paper in the School on Representation theory, Vertex and Chiral Algebras,
IMPA, Rio de Janeiro, March 14 - 18, 2022. The author would like to 
thank the organizers, especially, Jethro van Ekeren, for the invitation
and arrangement.

\renewcommand{\theequation}{\thesection.\arabic{equation}}
\renewcommand{\thethm}{\thesection.\arabic{thm}}
\setcounter{equation}{0} \setcounter{thm}{0}

\section{Rational functions, their generalizations and algebraic convergence}

The convergence of an expansion of a rational functions (or suitable simple 
generalizations of 
rational functions) is the simplest type of convergence appearing in the 
study of vertex operator algebras, their modules and 
twisted modules. This type of convergence is simple, algebraic and 
very useful, but unfortunately 
does not work for general intertwining operators. We begin our discussion of 
the convergence in this paper with this type of convergence 
and the reader should compare it with
the convergence to be discussed in later sections. 

Let 
$\C((z_{1}^{-1}, z_{2}))$ the space of all Laurent series of the form
$$\sum_{m, n\in \N}c_{mn}z_{1}^{-m_{0}-m}z_{2}^{n_{0}+n}$$
for $m_{0}, n_{0}\in \Z$ and $c_{mn}\in \C$ for $m, n\in \N$. 
Let $h(z_{1}, z_{2})\in \C((z_{1}^{-1}, z_{2}))$.
Assume that there is a nonnegative integer $N$ such that 
$$(z_{1}-z_{2})^{N}h(z_{1}, z_{2})
\in \C[z_{1}, z_{1}^{-1}, z_{2},  z_{2}^{-1}].$$
Then  we see that  $h(z_{1}, z_{2})$ 
is equal to the product of an element of 
$\C[z_{1}, z_{1}^{-1}, z_{2}, z_{2}^{-1}]$
 and  the expansion of $(z_{1}-z_{2})^{-N}$ in $\C((z_{1}^{-1}, z_{2}))$
 (that is, in nonnegative powers of the second variable $z_{2}$). 
 This product can always be written as the expansion 
 in $\C((z_{1}^{-1}, z_{2}))$ of 
 a rational function
 \begin{equation}\label{rat-fn}
 f(z_{1}, z_{2})=\frac{g(z_{1}, z_{2})}{z_{1}^{m}z_{2}^{n}(z_{1}-z_{2})^{l}},
 \end{equation}
 where $g(z_{1}, z_{2})\in \C[z_{1}, z_{2}]$ and $m, n, l\in \N$. 
 In other words,
 $h(z_{1}, z_{2})$ is absolutely convergent in the region $|z_{1}|>|z_{2}|>0$
 to the rational function $f(z_{1}, z_{2})$. The rational function $f(z_{1}, z_{2})$
 is an analytic function defined on the region $M^{2}=\{(z_{1}, z_{2})\in
 \C^{2}\mid z_{1}, z_{2}\ne 0, z_{1}\ne z_{2}\}$ and thus is 
 the analytic extension 
to $M^{2}$ of the sum of $h(z_{1}, z_{2})$ 
on the region $|z_{1}|>|z_{2}|>0$. 
 
 This simple fact gives a method to prove the convergence of a suitable series
 $h(z_{1}, z_{2})\in \C((z_{1}^{-1}, z_{2}))$:
 To prove this convergence, we need only find a
 nonnegative integer $N$ such that $(z_{1}-z_{2})^{N}h(z_{1}, z_{2})$
 is a Laurent polynomial in $z_{1}$ and $z_{2}$. 
 Note that the algebraic formulation of 
 this convergence holds  when we replace $\C$ by any field 
 $\F$ of characteristic $0$. 
 
 This method has been used extensively in the study of vertex operator algebras, 
 modules, twisted modules and some of their generalizations. 
 See for example the books
 \cite{DL} and \cite{LL} and the references there for details. 
 Here we use the product of two formal Laurent 
 series of operators acting on a lower-bounded graded vector space
 $V=\coprod_{n\in n_{0}+\N}V_{(n)}$
 to demonstrate the use of this method. 
 Let 
 $$A_{1}(z)=\sum_{n\in \Z}A_{1}^{n}z^{-n-1}, A_{2}(z)=\sum_{n\in \Z}A_{2}^{n}z^{-n-1}
 \in (\text{End}\;V)[[z, z_{1}]],$$
 and assume that $A_{1}^{n}$ and $A_{2}^{n}$ for $n\in \Z$ are operators on $V$ of weight (degree)
 $\lambda_{1}-n-1$ and $\lambda_{2}-n-1$, where $\lambda_{1}, \lambda_{2}\in \C$.
Then for $v\in V$ and $v'\in V'$,
 $\langle v', A_{1}(z)v\rangle, \langle v', A_{1}(z)v\rangle\in \C[z, z^{-1}]$, 
 $\langle v', A_{1}(z_{1})A_{2}(z_{2})v\rangle\in \C((z_{1}^{-1}, z_{2}))$ and 
 $\langle v', A_{2}(z_{2})A_{1}(z_{1})v\rangle\in \C((z_{2}^{-1}, z_{1}))$. 
 Assume that $A_{1}(z)$ and $A_{2}(z)$ satisfy the weak commutativity 
 (or locality), that is, 
 there exists $N\in \N$ such that 
\begin{equation}\label{weak-comm}
(z_{1}-z_{2})^{N}\langle v', A_{1}(z_{1})A_{2}(z_{2})v\rangle
=(z_{1}-z_{2})^{N}\langle v', A_{2}(z_{2})A_{1}(z_{1})v\rangle
\end{equation}
for $v\in V$ and $v'\in V'$.
Since $(z_{1}-z_{2})^{N}$ is a polynomial in $z_{1}$ and $z_{2}$, 
the left- and right-hand sides of (\ref{weak-comm}) are in 
$\C((z_{1}^{-1}, z_{2}))$ and $\C((z_{2}^{-1}, z_{1}))$, respectively. 
 Then the equality (\ref{weak-comm}) implies that both sides are in 
$\C[z_{1}, z_{1}^{-1}, z_{2},  z_{2}^{-1}]$. Using the method 
we discussed above, we see that 
$\langle v', A_{1}(z_{1})A_{2}(z_{2})v\rangle$ 
and $\langle v', A_{2}(z_{2})A_{1}(z_{1})v\rangle$
are absolutely convergent in the region $|z_{1}|>|z_{2}|>0$
and $|z_{2}|>|z_{1}|>0$, respectively,  to a common rational function
of the form (\ref{rat-fn}). 

We note that the application of this method depends heavily on the weak commutativity.
 
This method can be generalized immediately  to the case that 
$h(z_{1}, z_{2})$ is in 
\begin{equation}\label{ring}
\sum_{k_{1}, k_{2}=0}^{K}\sum_{j=1}^{J}z_{1}^{r_{j}}z_{2}^{s_{j}}
(\log z_{1})^{k_{1}}(\log z_{2})^{k_{2}}\C((z_{1}^{-1}, z_{2}))
\subset \C\{z_{1}, z_{2}\}[\log z_{1}, \log z_{2}],
\end{equation}
where $r_{j}, s_{j}\in \C$ for $j=1, \dots, J$. 
Assume that there is a nonnegative integer $N$ such that 
$$(z_{1}-z_{2})^{N}h(z_{1}, z_{2})
\in \sum_{k_{1}, k_{2}=0}^{K}
\sum_{j=1}^{J}z_{1}^{r_{j}}z_{2}^{s_{j}}
(\log z_{1})^{k_{1}}(\log z_{2})^{k_{2}}\C[z_{1}, z_{1}^{-1}, z_{2},  z_{2}^{-1}].$$
Then the same argument as above shows that 
$h(z_{1}, z_{2})$ 
is equal to the expansion  in (\ref{ring})
of  a function  
\begin{equation}\label{g-rat-fn}
f(z_{1}, z_{2})=\sum_{k_{1}, k_{2}=0}^{K}\sum_{j=1}^{J}
 \frac{g(z_{1}, z_{2})}{z_{1}^{m-r_{j}}z_{2}^{n-s_{j}}(z_{1}-z_{2})^{l}}
 (\log z_{1})^{k_{1}}(\log z_{2})^{k_{2}},
 \end{equation}
 where $g(z_{1}, z_{2})\in \C[z_{1}, z_{2}]$ and $m, n, l\in \N$. 
 The discussion above treats the series $h(z_{1}, z_{2})$ as a formal series
 and the function  $f(z_{1}, z_{2})$  as an element of the  localization of 
 the ring (\ref{ring}) by positive powers of $z_{1}-z_{2}$. If we use complex variables, we need to take 
 values of $z_{1}^{r_{j}}$, $z_{2}^{s_{2}}$, $\log z_{1}$ and $\log z_{2}$. 
 For any such values,  we see that $h(z_{1}, z_{2})$ evaluated at $z_{1}$ and 
 $z_{2}$ using these values is 
 absolutely convergent in the region $|z_{1}|>|z_{2}|>0$
 to the function (\ref{g-rat-fn}) evaluated using the same values of 
 $z_{1}^{r_{j}}$, $z_{2}^{s_{2}}$, $\log z_{1}$ and $\log z_{2}$. 
This generalization has been used in the study of twisted modules for 
a vertex operator algebras 
and also in the study of the product of one vertex operator for 
a module and one intertwining operator. 

There is also a generalization to 
abelian intertwining operator algebras introduced by Dong and Lepowsky 
in \cite{DL}. We still consider
$h(z_{1}, z_{2})$ in (\ref{ring}). But we assume that 
there is a complex number $N$ instead of a nonnegative integer such that 
$$(z_{1}-z_{2})^{N}h(z_{1}, z_{2})
\in \sum_{k_{1}, k_{2}=0}^{K}
\sum_{j=1}^{J}z_{1}^{r_{j}}z_{2}^{s_{j}}
(\log z_{1})^{k_{1}}(\log z_{2})^{k_{2}}\C[z_{1}, z_{1}^{-1}, z_{2},  z_{2}^{-1}],$$
where $(z_{1}-z_{2})^{N}$ is understood as the binomial expansion 
of $(z_{1}-z_{2})^{N}$ in nonnegative powers of $z_{2}$, that is, the 
expansion of $(z_{1}-z_{2})^{N}$ in the region $|z_{1}|>|z_{2}|>0$. 
Then the same argument as above shows that 
$h(z_{1}, z_{2})$ 
is equal to the expansion  as a series in powers of $z_{1}$ and $z_{2}$ 
with only finitely many negative real parts of the powers of $z_{2}$
of  a function  
$$f(z_{1}, z_{2})=\sum_{k_{1}, k_{2}=0}^{K}\sum_{j=1}^{J}
 \frac{g(z_{1}, z_{2})}{z_{1}^{m-r_{j}}z_{2}^{n-s_{j}}(z_{1}-z_{2})^{l}}
 (\log z_{1})^{k_{1}}(\log z_{2})^{k_{2}},$$
 where $g(z_{1}, z_{2})\in \C[z_{1}, z_{2}]$, $m, n\in \N$ and $l\in \C$. 
 In other words, $h(z_{1}, z_{2})$ evaluated using any values
 $z_{1}^{r_{j}}$, $z_{2}^{s_{2}}$, $(z_{1}-z_{2})^{-l}$,
 $\log z_{1}$ and $\log z_{2}$ is 
 absolutely convergent in the region $|z_{1}|>|z_{2}|>0$
 to the function $f(z_{1}, z_{2})$ evaluated using the same values. 
 But in practice, it is not easy to apply this method in this general case
 since the generalization of the weak commutativity in this case is not 
easy to verify. In fact, since the power $N$ is not a nonnegative integer 
anymore, the expansion of 
$(z_{1}-z_{2})^{N}$ in the regions $|z_{1}|>|z_{2}|>0$
and $|z_{2}|>|z_{1}|>0$ are very different. In this case,
the weak commutativity for two (multivalued) fields 
$A_{1}(z_{1})$ and $A_{2}(z_{2})$
is of the form
$$(z_{1}-z_{2})^{N}A_{1}(z_{1})A_{2}(z_{2})
=(-z_{2}+z_{1})^{N}A_{2}(z_{2})A_{1}(z_{1}),$$
where $(z_{1}-z_{2})^{N}$ is the binomial expansion 
of $(z_{1}-z_{2})^{N}$ in nonnegative powers of $z_{2}$ and 
$(-z_{2}+z_{1})^{N}$ is  the binomial expansion 
of $(-z_{2}+z_{1})^{N}$ in nonnegative powers of $z_{1}$.
In general, this weak commutativity is not always easy to verify since 
$(z_{1}-z_{2})^{N}$ and $(-z_{2}+z_{1})^{N}$ are very different.

For general intertwining operators (for example, intertwining 
operators among modules for the affine vertex operator algebras 
and Virasoro vertex operator algebras), 
even weak commutativity for abelian 
intertwining operator algebras is not satisfied.
This is the main reason why we need the method in the next section to 
prove the convergence of products of intertwining operators. 
Even in some special cases that we expect the intertwining operators 
to form an abelian intertwining operator algebra, 
because of the difficulty to verify the 
weak commutativity mentioned above, we still need the method 
in the next section to prove the the convergence of 
products of intertwining operators.

 Our discussions above are for series in only two complex variables 
 or products of two series of operators. But they
 can be generalized easily to the case of an arbitrary number of complex 
 variables or products of an arbitrary number of series of operators.

\renewcommand{\theequation}{\thesection.\arabic{equation}}
\renewcommand{\thethm}{\thesection.\arabic{thm}}
\setcounter{equation}{0} \setcounter{thm}{0}

\section{Convergence of products of 
intertwining operators}

 The convergence discussed in the preceding section does not work for 
 general intertwining operators. We have to use a completely different 
 method. In this section, 
 we discuss the convergence of products of intertwining operators
 and their analytic extensions using this method.
 
Let $V$ be a vertex operator algebra, 
$W_{1}$, $W_{2}$, $W_{3}$, $W_{4}$ and $W_{5}$ 
 be generalized $V$-modules and $\Y_{1}$ and $\Y_{2}$ (logarithmic)
 intertwining operators of type $\binom{W_{4}}{W_{1}W_{5}}$ 
 and $\binom{W_{5}}{W_{2}W_{3}}$,
 respectively. For $w_{1}\in W_{1}$ and $w_{2}\in W_{2}$, 
 $$\Y_{1}(w_{1}, z_{1})\in \hom(W_{5}, W_{4})[\log z_{1}]\{z_{1}\}$$ 
 and 
$$\Y_{1}(w_{2}, z_{2})\in \hom(W_{3}, W_{5})[\log z_{2}]\{z_{2}\}.$$
 Then 
 $$\Y_{1}(w_{1}, z_{1})
 \Y_{1}(w_{2}, z_{2})\in \hom(W_{3}, W_{4})[\log z_{2}, \log z_{2}]
 \{z_{1}, z_{2}\}
 .$$
 The first problem in the study of the product $\Y_{1}(w_{1}, z_{1})
 \Y_{1}(w_{2}, z_{2})$ is the convergence of this product in a suitable sense. 
 Since the series $\Y_{1}(w_{1}, z_{1})
 \Y_{1}(w_{2}, z_{2})$ contains nonintegral powers of $z_{1}$ and $z_{2}$ and 
nonnegative integral powers of 
$\log z_{1}$ and $\log z_{2}$, we first need to choose values of 
 these powers of $z_{1}$ and $z_{2}$ and values of $\log z_{1}$ and $\log z_{2}$.
In fact, if the values of $\log z_{1}$ and $\log z_{2}$ are chosen to be 
$l_{p}(z_{1})=\log |z_{1}|+i \arg z_{1}+2\pi i p$ and 
$l_{q}(z_{2})=\log |z_{2}|+i \arg z_{2}+2\pi i q$, where 
$0 \le \arg z_{1}, \arg z_{2}<2\pi$,  then these values also
give values $e^{ml_{p}(z_{1})}$ and $e^{nl_{q}(z_{2})}$
of the powers $z_{1}^{m}$ and $z_{2}^{n}$ 
of $z_{1}$ and $z_{2}$. Using these values, we obtain a series in $\C$
 \begin{equation}\label{int-op-prod}
 \langle w_{4}', \Y_{1}(w_{1}, z_{1})
 \Y_{1}(w_{2}, z_{2})w_{3}\rangle\lbar_{\log z_{1}=l_{p}(z_{1}),\; \log z_{2}=l_{p}(z_{2})}
 \end{equation}
 for $w_{1}\in W_{1}$, $w_{2}\in W_{2}$, 
 $w_{3}\in W_{3}$ and $w_{4}'\in W_{4}'$. We want to know whether
 (\ref{int-op-prod}) 
 is absolutely convergent in a suitable region for $z_{1}$ and $z_{2}$.
 
 In general, (\ref{int-op-prod}) might not be convergent. 
 The method used in the preceding section 
 works only when it is possible to multiply a nonnegative integral powers 
 $(z_{1}-z_{2})^{N}$ of $z_{1}-z_{2}$ 
 to get a finite sum. But in general (even in the case of abelian intertwining 
 operator algebras), there is no such $N$; in general, 
 there is even no polynomial
 in $z_{1}$ and $z_{2}$ that can be multiplied 
 to (\ref{int-op-prod})
 to get a finite sum.

The method used to prove the convergence of (\ref{int-op-prod})
is to show that the series (\ref{int-op-prod}) satisfies the expansion 
in the region $|z_{1}|>|z_{2}|>0$ 
of a system of differential equations
with coefficients in the ring 
$$R=\C[z_{1}, z_{1}^{-1}, z_{2}, z_{2}^{-1},
(z_{1}-z_{2})^{-1}]$$
and with regular singular points at $(z_{1}, z_{2})=(\infty, 0)$. 
In fact, since each coefficient of  (\ref{int-op-prod}) as a series in 
powers of $z_{2}$ and $\log z_{2}$ 
is a finite sum, we need only prove the convergence of 
(\ref{int-op-prod}) for fixed $z_{1}\in \C^{\times}$. In particular, we need only 
derive a differential equation in the variable $z_{2}$ 
with the regular singular point $z_{2}=0$. 
Then by the theory of differential equations of regular singular points, 
the formal series solution of the system of the differential equations of regular
singular points must be the expansion of an analytic solution of the system. 
In other words, the series (\ref{int-op-prod}) is absolutely convergent 
in the region $|z_{1}|>|z_{2}|>0$ to an analytic solution. Since the coefficients 
of the system of differential equations are in the ring $R$, 
this solution on the region $|z_{1}|>|z_{2}|>0$ 
can be analytically extended to a multivalued analytic function on $M^{2}$. 
Moreover, this multivalued analytic extension also satisfies a system of differential equations
with coefficients in $R$ and with regular singular points at $(z_{2}, z_{1}-z_{2})=(\infty, 0)$. 
Then the multivalued analytic extension can be expanded 
in the region $|z_{1}|>|z_{1}-z_{2}|>0$ as a series containing 
terms in powers of $z_{2}$ and $z_{1}-z_{2}$ and nonnegative 
integral powers of logarithms of 
$z_{1}$ and $z_{2}$.

 In the special case of Wess-Zumino-Witten models or minimal models, 
 we have the Knizhnik-Zamolodchikov equations
 \cite{KZ}  or the Belevin-Polyakov-Zamolodchikov equations \cite{BPZ},
 respectively.
 The Knizhnik-Zamolodchikov equations were used by Tsuchiya and Kanie 
 \cite{TK} to 
 prove the convergence of products of intertwining operators 
 (called vertex operators in \cite{TK}) among suitable modules for 
 the affine Lie algebra $\widehat{\mathfrak{sl}(2, \C)}$.
 The  Belevin-Polyakov-Zamolodchikov equations
 and Knizhnik-Zamolodchikov equations
 were used by the author \cite{H-minimal-mod} and 
 by Lepowsky and the author \cite{HL-affine}, respectively, 
 to prove the convergence of products 
 of intertwining  operators for vertex operator algebras
 for the minimal models and for the Wess-Zumino-Witten models. 
 But for general 
 vertex operator algebras, products of intertwining operators
 might not satisfy differential equations. Some conditions on 
 the vertex operator algebras and modules must be satisfied 
 in order to have such differential equations.

Note that  for a solution of a system of differential equations with coefficients 
in $R$, the derivatives of the solution must span
a finitely-generated module over the ring $R$. In particular, if (\ref{int-op-prod}) indeed 
converges absolutely to a solution of such a system of differential equations,
the derivatives of (\ref{int-op-prod}) span a finitely-generated module over 
$R$. Together with the $L(-1)$-derivative property of intertwining operators, 
it is not difficult to see from this fact that there should be some finiteness conditions satisfied by 
the grading-restricted generalized $V$-modules $W_{1}$, $W_{2}$, $W_{3}$ and $W_{4}'$. 
This is the reason why in \cite{H-diff-eqn}, the following $C_{1}$-cofiniteness condition
on $W_{1}$, $W_{2}$, $W_{3}$ and $W_{4}'$ is needed:
For a grading-restricted generalized $V$-module $W$, 
we say that $W$ is $C_{1}$-cofinite if 
 $\dim W/C_{1}(W)<\infty$, where $C_{1}(V)$ is the subspace of $W$ spanned by 
 elements of the form $\res_{x}x^{-1}Y_{W}(v, x)w$ for 
 $v\in V_{+}=\coprod_{n\in \Z_{+}}V_{(n)}$ and $w\in W$. We also need 
 another condition on generalized $V$-modules
  in our precise statement of the theorem below:
  A generalized $V$-module is said to be quasi-finite dimensional if
  for any $N\in \R$, the subspace $\coprod_{\Re(n)\le N}W_{[n]}$
  is finite dimensional. In the case that all irreducible generalized $V$-modules 
  are grading restricted, a generalized $V$-modules of finite length 
  must be quasi-finite dimensional.

 Using the the Jacobi identity for
 intertwining operators and vertex operators acting on modules, 
 we obtain certain identities for series of the form 
 (\ref{int-op-prod}). For example, we have 
 \begin{align*}
&\langle w_{4}', \Y_1(w_1, z_1) 
\Y_2(\res_{x}x^{-1}Y_{W_{1}}(u, x)w_2, z_2)w_{3}\rangle\nn
&\quad =\sum_{k\in \N}z_{2}^{k} 
\langle (\res_{x}x^{-1-k}Y_{W_{4}}(u, x))'w_{4}', 
\Y_1(w_1, z_1) \Y_2(w_2, z_2)w_{3}\rangle\nn
& \quad\quad +\sum_{k\in \N}(-1)^{k}
(z_{1}-z_{2})^{-1-k} \langle w_{4}', 
\Y_1(\res_{x}x^{k}Y_{W_{2}}(u, x)w_1, z_1) 
\Y_2(w_2, z_2)w_{3}\rangle \nn
&\quad\quad +\sum_{k\in \N}z_{2}^{-1-k}
\langle w_{4}', \Y_1(w_1, z_1) \Y_2(w_2, z_2)
\res_{x}x^{k}Y_{W_{2}}(u, x)w_{3}\rangle,
\end{align*}
where $(\res_{x}x^{-1-k}Y_{W_{2}}(u, x))'$ is the adjoint on 
$W_{4}'$ of 
$\res_{x}x^{-1-k}Y_{W_{4}}(u, x)$. 
 Using these identities, the $C_{1}$-cofiniteess condition, the 
 quasi-finite-dimensionality of the generalized $V$-modules involved
 and the $L(-1)$-derivative property,
 it was proved in \cite{H-diff-eqn} that the derivatives 
 of (\ref{int-op-prod}) span a finitely-generated module over 
 $R$ and thus must satisfy the expansion in the region $|z_{1}|>|z_{2}|>0$ 
of a system of differential equations
with coefficients in $R$. Using more careful examinations of the coefficients of 
the differential equations, it was proved in \cite{H-diff-eqn} that  one can 
always  find such a system of differential equations such that 
the singular point $(z_{1}, z_{2})=(\infty, 0)$ or $(z_{2}, z_{1}-z_{2})=(\infty, 0)$
is regular. (As is mentioned above, in fact it is enough to show that 
for fixed $z_{1}\in \C^{\times}$, (\ref{int-op-prod}) satisfies 
a differential equation in the variable $z_{2}$ 
with the regular singular point $z_{2}=0$.)
Then we obtain the convergence and analytic extension result
for products of two intertwining operators. 

 Below is the precise statement of the convergence and analytic extension result.
It is essentially the $n=2$ case of Theorem 11.8 in \cite{HLZ7} 
with the category $\mathcal{C}$
 being the category of grading-restricted generalized $V$-modules. Its 
 proof was in fact given in the proof of Theorem 3.5 in \cite{H-diff-eqn}, where
 it is proved in addition that when all $\N$-gradable weak $V$-modules are completely 
 reducible, there is no logarithms of the variables in the expansion near the 
 singular point $(z_{2}, z_{1}-z_{2})=(\infty, 0)$.
 
 \begin{thm}[\cite{H-diff-eqn}, \cite{HLZ7}]\label{conv-prod-int}
 Let $V$ be a vertex operator algebra, $W_{1}$, $W_{2}$, $W_{3}$, $W_{4}$ and $W_{5}$ 
 be generalized $V$-modules and $\Y_{1}$ and $\Y_{2}$ (logarithmic)
 intertwining operators of type $\binom{W_{4}}{W_{1}W_{5}}$ 
 and $\binom{W_{5}}{W_{2}W_{3}}$,
 respectively. Assume that $W_{1}$, $W_{2}$, $W_{3}$, $W_{4}'$ are 
 quasi-finite dimensional and $C_{1}$-cofinite. Then
 (\ref{int-op-prod})
is absolutely convergent when $|z_1|>|z_2|>0$ and can be analytically
extended to a multivalued analytic function on $M^{2}$.
\end{thm}

The convergence and analytic extensions of products of more than two intertwining operators 
are proved similarly.  The convergence and analytic extensions of iterates
of intertwining operators can be derived from Theorem \ref{conv-prod-int}
and properties of intertwining operators.
See \cite{H-diff-eqn}, \cite{HLZ6} and \cite{HLZ7}.

\renewcommand{\theequation}{\thesection.\arabic{equation}}
\renewcommand{\thethm}{\thesection.\arabic{thm}}
\setcounter{equation}{0} \setcounter{thm}{0}

\section{Convergence of $q$-traces and pseudo-$q$-traces of 
products of intertwining operators}

The multivalued analytic functions obtained from the analytic extensions of 
products of intertwining operators discussed in the preceding section are 
in fact the  genus-zero correlation functions for the corresponding 
chiral conformal field theories. 
To construct genus-one correlation functions from these 
genus-zero correlation functions, we need to take $q$-traces and 
pseudo-$q$-traces of 
products of intertwining operators. Since grading-restricted generalized $V$-modules 
in general are always infinite dimensional, the first problem one needs to solve 
is the convergence of these $q$-traces and pseudo-$q$-traces. 

As in the case of products of intertwining operators, we shall discuss only 
the case of $q$-traces and pseudo-$q$-traces of products of two intertwining operators. 
The general case is similar. 

Geometrically, products of intertwining operators correspond to 
genus-zero Riemann surfaces with punctures and local coordinates
 (see \cite{H-geom-voa} and  \cite{H-geom-int}) 
 but the $q$-traces or pseudo-$q$-traces
of products of intertwining operators correspond  to 
genus-one surfaces with punctures and local coordinates. 
Since the standard description of genus-one Riemann surfaces is in terms 
of parallelograms in the complex plane, not annuli in the sphere, 
to use intertwining operators to write down genus-one correlation functions, 
we have to modify intertwining operators correspondingly. 

For a grading-restricted generalized $V$-module $W$,
as in \cite{H-mod-inv-int}, let 
$$\mathcal{U}_{W}(x)=(2\pi ix)^{L_{W}(0)}e^{-L^{+}_{W}(A)}\in (\mbox{\rm End}\;W)\{x\}[\log x],$$
where $(2\pi i)^{L(0)}=e^{(\log 2\pi +i \frac{\pi}{2})L(0)}$, $x^{L_{W}(0)}
=x^{L_{W}(0)_{S}}e^{(\log x)L_{W}(0)_{N}}$ 
($L_{W}(0)_{S}$ and $L_{W}(0)_{N}$ being the semisimple and 
nilpotent, respectively,  parts of $L_{W}(0)$), 
$L^{+}_{W}(A)=\sum_{j\in \mathbb{Z}_{+}}
A_{j}L(j)$ and $A_{j}$ for $j\in \N$ are given by 
$$\frac{1}{2\pi i}\log(1+2\pi i y)=\left(\exp\left(\sum_{j\in \mathbb{Z}_{+}}
A_{j}y^{j+1}\frac{\partial}{\partial y}\right)\right)y.$$

Let 
$W_{1}, W_{2}, W_{3}, W_{4}$
be grading-restricted generalized $V$-modules.
Let $\Y_{1}$ and $\Y_{2}$ be intertwining operators of types
$\binom{W_{3}}{W_{1}W_{4}}$ and 
$\binom{W_{4}}{W_{2}W_{3}}$, respectively. For $z\in \C$, 
we shall use $q_{z}$ to denote $e^{2\pi iz}$. 
We call $\mathcal{Y}_{1}(\mathcal{U}_{W_{1}}(q_{z_{1}})w_{1}, 
q_{z_{1}})$ and 
$\mathcal{Y}_{2}(\mathcal{U}_{W_{2}}(q_{z_{2}})w_{2}, q_{z_{2}})$
for $w_{1}\in W_{1}$ and $w_{2}\in W_{2}$
geometrically-modified intertwining operators. 
For $z_{1}, z_{2}, \tau\in \C$, we have the $q$-trace (shifted by $-\frac{c}{24}$) 
\begin{align}\label{q-trace}
&\tr_{W_{3}}\mathcal{Y}_{1}(\mathcal{U}_{W_{1}}(q_{z_{1}})w_{1}, 
q_{z_{1}})
\mathcal{Y}_{2}(\mathcal{U}_{W_{2}}(q_{z_{2}})w_{2}, q_{z_{2}})
q^{L(0)-\frac{c}{24}}\nn
&\quad =\sum_{n\in \C}\tr_{(W_{3})_{[n]}}\left(\pi_{n}\mathcal{Y}_{1}(\mathcal{U}_{W_{1}}(q_{z_{1}})w_{1}, 
q_{z_{1}})
\mathcal{Y}_{2}(\mathcal{U}_{W_{2}}(q_{z_{2}})w_{2}, q_{z_{2}})
q^{L(0)-\frac{c}{24}}\lbar_{(W_{3})_{[n]}}\right)\nn
&\quad =\sum_{n\in \C}\tr_{(W_{3})_{[n]}}\left(\pi_{n}\mathcal{Y}_{1}(\mathcal{U}_{W_{1}}(q_{z_{1}})w_{1}, 
q_{z_{1}})
\mathcal{Y}_{2}(\mathcal{U}_{W_{2}}(q_{z_{2}})w_{2}, q_{z_{2}})
q^{n-\frac{c}{24}}e^{(\log q)L_{W_{3}}(0)_{N}}
\lbar_{(W_{3})_{[n]}}\right),
\end{align}
where $\pi_{n}: W_{3}\to (W_{3})_{[n]}$ for $n\in \C$ is the projection from 
the algebraic completion 
$$\overline{W}_{3}=\prod_{n\in \C}(W_{3})_{[n]}$$ of 
$$W_{3}=\coprod_{n\in \C}(W_{3})_{[n]}$$ 
to $(W_{3})_{[n]}$.

In general, we need to consider pseudo-$q$-traces of 
$$\mathcal{Y}_{1}(\mathcal{U}_{W_{1}}(q_{z_{1}})w_{1}, 
q_{z_{1}})
\mathcal{Y}_{2}(\mathcal{U}_{W_{2}}(q_{z_{2}})w_{2}, q_{z_{2}}).$$
For a pseudo-$q$-trace, we need to consider a 
grading-restricted generalized $V$-module
equipped with a projective right module structure for 
a finite-dimensional associative algebra $P$ over $\C$. 
We first define the pseudo-trace of an operator $\alpha\in \text{End}_{P}M$ 
on a finitely generated projective right $P$-module $M$. 
Since $P$ is projective, 
for such a right $P$-module $M$, there exists a projective basis, that is,
a pair of sets $\{w_i\}_{i=1}^n\subset M$,
$\{w'_i\}_{i=1}^n\subset \hom_P(M, P)$
such that for all $w\in M$, $w = \sum_{i=1}^n w_i(w'_i(w)).$
A linear function $\phi: P\rightarrow \C$ is said to be symmetric
if  $\phi(pq) = \phi(qp)$ for all $p,q\in P$.
For a symmetric linear function $\phi$,
the pseudo-trace  $\tr^{\phi}_M\alpha$ of 
$\alpha\in \text{End}_{P}M$ associated to $\phi$ 
is defined by
$$\tr^{\phi}_M \alpha = \phi\left(\sum_{i=1}^n w'_i(\alpha(w_i))\right).$$ 
For a grading-restricted generalized $V$-module $W$ 
equipped with a projective right $P$-module structure, its homogeneous subspaces $W_{[n]}$ 
for $n\in \C$ are finitely generated projective right $P$-modules. Then for a given 
symmetric linear function $\phi$ on $P$, we have the pseudo-trace $\tr^{\phi}_M\alpha_{n}$
of $\alpha_{n}\in {\rm End}_{P}W_{[n]}$. For $\alpha\in {\rm End}_{P}W$, we define 
the pseudo-$q$-trace  (shifted by $-\frac{c}{24}$)
of $\alpha$ by 
\begin{align*}
\tr^{\phi}_W\alpha q^{L_{W}(0)-\frac{c}{24}}
&=\sum_{n\in \C}\tr^{\phi}_{W_{[n]}}
\left(\pi_{n}\alpha
q^{L_{W}(0)-\frac{c}{24}}\lbar_{W_{[n]}}\right)\nn
&=\sum_{n\in \C}\tr^{\phi}_{W_{[n]}}
\left(\pi_{n}\alpha
q^{n-\frac{c}{24}}e^{L_{W}(0)_{N}\log q}\lbar_{W_{[n]}}\right).
\end{align*}
Note that $\tr^{\phi}_W\alpha q^{L_{W}(0)-\frac{c}{24}}$ 
defined above is a series. To obtain the pseudo-$q$-trace of $\alpha$
as a function of $q$, we have to prove its convergence. 

What we are interested is the pseudo-$q$-trace  (shifted by $-\frac{c}{24}$)
\begin{align}\label{pseudo-q-trace}
&\tr^{\phi}_{W_{3}}\mathcal{Y}_{1}(\mathcal{U}_{W_{1}}(q_{z_{1}})w_{1}, 
q_{z_{1}})
\mathcal{Y}_{2}(\mathcal{U}_{W_{2}}(q_{z_{2}})w_{2}, q_{z_{2}})
q^{L(0)-\frac{c}{24}}\nn
&\quad =\sum_{n\in \C}\tr^{\phi}_{(W_{3})_{[n]}}
\left(\pi_{n}\mathcal{Y}_{1}(\mathcal{U}_{W_{1}}(q_{z_{1}})w_{1}, 
q_{z_{1}})
\mathcal{Y}_{2}(\mathcal{U}_{W_{2}}(q_{z_{2}})w_{2}, q_{z_{2}})
q^{L_{W_{3}}(0)-\frac{c}{24}}\lbar_{(W_{3})_{[n]}}
\right)\nn
&\quad =\sum_{n\in \C}\tr^{\phi}_{(W_{3})_{[n]}}
\left(\pi_{n}\mathcal{Y}_{1}(\mathcal{U}_{W_{1}}(q_{z_{1}})w_{1}, 
q_{z_{1}})
\mathcal{Y}_{2}(\mathcal{U}_{W_{2}}(q_{z_{2}})w_{2}, q_{z_{2}})
q^{n-\frac{c}{24}}e^{(\log q) L_{W_{3}}(0)_{N}}\lbar_{(W_{3})_{[n]}}
\right),
\end{align}
where $\phi$ is a symmetric linear function
on a finite-dimensional associative algebra $P$, $W_{3}$
is a projective right $P$-module such that its vertex operators and
$$\mathcal{Y}_{1}(\mathcal{U}_{W_{1}}(q_{z_{1}})w_{1}, 
q_{z_{1}})
\mathcal{Y}_{2}(\mathcal{U}_{W_{2}}(q_{z_{2}})w_{2}, q_{z_{2}})$$
commute with the action of $P$ on $W_{3}$.
Note that the $q$-trace (\ref{q-trace}) is the special case 
of (\ref{pseudo-q-trace}) for which $P=\{e\}$, where $e$ is the identity
element of $P$, and $\phi$ is given by $\phi(e)=1$. 
We want to know whether (\ref{q-trace}) and (\ref{pseudo-q-trace}) 
are absolutely convergent in a suitable region for $z_{1}, z_{2}$ and $q$. 

In general, (\ref{q-trace}) and (\ref{pseudo-q-trace}) might not be convergent. 
Just as in the case of 
 products of intertwining operators,
we also need some cofiniteness condition. 
When the confiniteness condition is satisfied,
 (\ref{q-trace}) and (\ref{pseudo-q-trace})  are convergent in a suitable region.
 The method that we use to prove this convergence is to show that 
 the series (\ref{q-trace}) and (\ref{pseudo-q-trace}) satisfy the expansion 
 in the region $1>|q_{z_{1}}|>|q_{z_{2}}|>|q|>0$
 of a system of differential equations with coefficients in the ring 
 $\C[G_{4}(\tau), G_{6}(\tau), \wp_{2}(z_{1}-z_{2}; \tau), 
 \wp_{3}(z_{1}-z_{2}; \tau)]$,
 where 
 \begin{align*}
 G_{4}(\tau)&=\sum_{(k, l)\ne (0, 0)}\frac{1}{(k\tau+l)^{4}}, \\
 G_{6}(\tau)&=\sum_{(k, l)\ne (0, 0)}\frac{1}{(k\tau+l)^{6}}
 \end{align*}
 are Eisenstein series and 
 \begin{align*}
\wp_{2}(z; \tau)&=\frac{1}{z}+\sum_{(k, l)\ne (0, 0)}\left(\frac{1}{(z-(k\tau+l))^{2}}
-\frac{1}{(k\tau+l)^{2}}\right), \\
 \wp_{3}(z; \tau)&=-\frac{1}{2}\frac{\partial}{\partial z}\wp_{2}(z; \tau)
 \end{align*}
 are Weierstrass $\wp$-function and its derivative with respect to $z$ 
 multiplied by $-\frac{1}{2}$. 

The first convergence result on $q$-traces of products of vertex operators 
on a $V$-module was obtained by Zhu \cite{Z} for $V$ satisfying the
conditions  that $V$ has no 
nonzero element of negative weights, every lower-bounded generalized $V$-module
is completely reducible and a cofiniteness condition which is now called 
the $C_{2}$-cofiniteness condition. 
Let $C_{2}(V)$ be the subspace of $V$ spanned by elements of the form 
$\res_{x}x^{-2}Y_{V}(u, v)v$ for $u. v\in V$. If $\dim V/C_{2}(V)<\infty$, 
we say that $V$ is $C_{2}$-cofinite. 
The convergence result of Zhu was generalized by Miyamoto in \cite{M1} 
to a convergence result on 
$q$-traces of products of  one intertwining operator
and several vertex operators on modules in the case that $V$ satisfies the three conditions mentioned above 
in Zhu's paper \cite{Z}. It was also generalized 
by Miyamoto in \cite{M}
to a convergence result on pseudo-$q$-traces of products of vertex operators 
on a $V$-module in the case that $V$ has no 
nonzero element of negative weights and is $C_{2}$-cofinite. 
These convergence results are all 
proved in two steps: (i) An algebraic recurrence relation is
proved to reduce
the convergence of (pseudo-)$q$-traces of products of $n$ vertex operators 
on a $V$-module 
(or products of $n-1$ vertex operators on a $V$-module
and one intertwining operator)
to the convergence of (pseudo-)$q$-traces of one vertex operator
(or one intertwining operator)
on the same $V$-module. (ii) The convergence of 
(pseudo-)$q$-traces of one vertex operator (or one intertwining operator) 
on the same $V$-module is proved by using differential equations of 
regular singular points. 

The proofs of the algebraic recurrence relations in step (i) above 
need the 
commutator formula for vertex operators on $V$-modules
or the commutator formula 
between vertex operators on $V$-modules and intertwining operators. 
Since in general there is no commutator formula for intertwining operators,
there is no algebraic recurrence relation to reduce
the convergence of (\ref{q-trace}) and (\ref{pseudo-q-trace})
to the convergence of (pseudo-)$q$-traces of intertwining operators.
This difficulty is the main reason why the modular invariance of 
the space of $q$-traces of products of at least two intertwining operators
had been a conjecture for many years after Zhu's work \cite{Z}. 

In \cite{H-mod-inv-int}, the author proved the convergence of (\ref{q-trace}). 
In \cite{F1} and \cite{F2}, using the same method, 
Fiordalisi generalized the convergence of (\ref{q-trace})
proved in \cite{H-mod-inv-int} to the convergence of
(\ref{pseudo-q-trace}). In fact, 
certain formulas on $q$-traces of 
products of geometrically-modified intertwining operators 
are proved in \cite{H-mod-inv-int} and generalized 
in \cite{F1} and \cite{F2} to pseudo-$q$-traces of 
products of geometrically-modified intertwining operators.
For example, one such formula is 
\begin{align*}
&\tr^{\phi}_{W_{3}}
\mathcal{Y}_{1}(\mathcal{U}(q_{z_{1}})w_{1}, q_{z_{1}})
\mathcal{Y}_{2}(\mathcal{U}(q_{z_{2}})\res_{x}x^{-2}
Y_{W_{2}}(u, x)w_{2}, q_{z_{2}})q^{L(0)-\frac{c}{24}}\nn
&\quad=-\sum_{k\in \Z_{+}}(2k+1)\tilde{G}_{2k+2}(q)
\tr^{\phi}_{W_{3}}\mathcal{Y}_{1}(\mathcal{U}(q_{z_{1}})w_{1}, q_{z_{1}})
\mathcal{Y}_{2}(\mathcal{U}(q_{z_{2}})\res_{x}x^{2k}
Y_{W_{2}}(u, x)w_{2}, q_{z_{2}})q^{L(0)-\frac{c}{24}}\nn
&\quad\quad-\sum_{m\in \N}(-1)^{m}
(m+1)\tilde{\wp}_{m+2}(z_{i}-z_{j}; q)\cdot\nn
&\quad\quad\quad\quad\quad\quad\cdot 
\tr^{\phi}_{W_{3}}
\mathcal{Y}_{1}(\mathcal{U}(q_{z_{1}})\res_{x}x^{m}
Y_{W_{1}}(u, x)w_{1}, q_{z_{1}})
\mathcal{Y}_{2}(\mathcal{U}(q_{z_{2}})w_{n}, q_{z_{2}})
q^{L(0)-\frac{c}{24}},
\end{align*}
where $\tilde{G}_{2k+2}(q)$ for $k\in \Z_{+}$ 
are the $q$-expansions of the Eisenstein series $G_{2k+2}(\tau)$
and $\tilde{\wp}_{m+2}(z; q)$ for $m\in \N$
are the $q$-expansions of 
$$\wp_{m+2}(z; \tau)=\frac{(-1)^{m}}{(m+1)!}
\frac{\partial^{m}}{\partial z^{m}}\wp_{2}(z; \tau).$$
These formulas together with 
the $C_{2}$-cofiniteness of $V$ are used in
\cite{H-mod-inv-int}, \cite{F1} and \cite{F2} to show that the modules 
over  $\C[G_{4}(\tau), G_{6}(\tau), \wp_{2}(z_{1}-z_{2}; \tau), 
 \wp_{3}(z_{1}-z_{2}; \tau)]$
generated by these (pseudo-)$q$-traces 
are finitely generated. 
Another formula involving 
$q\frac{\partial}{\partial q}$, $\frac{\partial}{\partial z_{1}}$,
$\frac{\partial}{\partial z_{2}}$,
$L_{W_{1}}(0)$ and $L_{W_{2}}(0)$, 
the $q$-expansion of the Eisenstein series 
$G_{2}(\tau)$ and the $q$-expansion 
of the Weierstrass zeta-function $\wp_{1}(z; \tau)$ was
also proved in \cite{H-mod-inv-int} and 
generalized in \cite{F1} and \cite{F2}. We refer the reader to 
\cite{H-mod-inv-int}, \cite{F1} and \cite{F2} for this formula.
This formula  in fact gives how a modular invariant differential 
operator containing 
$q\frac{\partial}{\partial q}$, $\frac{\partial}{\partial z_{1}}$ and
$\frac{\partial}{\partial z_{2}}$ acts on $q$-traces and
pseudo-$q$-traces of 
products of geometrically-modified intertwining operators.
In \cite{H-mod-inv-int}, 
\cite{F1} and \cite{F2}, the action of this differential operator 
and the result that 
the module over $\C[G_{4}(\tau), G_{6}(\tau), \wp_{2}(z_{1}-z_{2}; \tau), \wp_{3}(z_{1}-z_{2}; \tau)]$ generated by these $q$-traces and
pseudo-$q$-traces
is finitely generated are used to prove that 
 (\ref{q-trace}) and (\ref{pseudo-q-trace})
satisfy the expansion in the region $1>|q_{z_{1}}|>|q_{z_{2}}|>|q_{\tau}|>0$
of a modular invariant system of differential equations with coefficients in 
$\C[G_{4}(\tau), G_{6}(\tau), \wp_{2}(z_{1}-z_{2}; \tau), \wp_{3}(z_{1}-z_{2}; \tau)]$. 
The system 
can also be chosen to be of regular singular point at each singular point. 
Then the convergence of (\ref{q-trace}) and (\ref{pseudo-q-trace}) follows. 

Below is the precise statement of the convergence and analytic extension 
of (\ref{pseudo-q-trace}). As is mentioned above,
 (\ref{q-trace}) is a special case of 
(\ref{pseudo-q-trace}). In the semisimple case that 
every lower-bounded generalized $V$-module is completely reducible,
(\ref{pseudo-q-trace}) is the same as (\ref{q-trace}). 

\begin{thm}[\cite{H-mod-inv-int}, \cite{F1}, \cite{F2}]
Assume that the vertex operator algebra $V$ has no 
nonzero element of negative weights and is $C_{2}$-cofinite. Then
in the region $1>|q_{z_{1}}|
>|q_{z_{2}}|>|q_{\tau}|>0$, the series
(\ref{pseudo-q-trace}) with $q=q_{\tau}=e^{2\pi i\tau}$
is absolutely convergent and can be analytically extended
to a multivalued analytic function in the region given by
$\Im(\tau)>0$ (here $\Im(\tau)$ is the imaginary part of $\tau$), 
$z_{1}\ne z_{2}+k\tau+l$ for $k, l\in \mathbb{Z}$. Moreover, the singular point
$z_{1}=z_{2}+k\tau+l$ for each $k, l\in \mathbb{Z}$ is regular, that is,
any branch of the multivalued analytic function can be expanded in 
a neighborhood of 
the singular point $z_{1}=z_{2}+k\tau+l$ as a series of the form
$$\sum_{p=0}^{K}
\sum_{j=1}^{M}(z_{1}-z_{2}+k\tau+l)^{r_{j}}(\log (z_{1}-z_{2}+k\tau+l))^{p}
f_{j, p}(z_{1}-z_{2}+k\tau+l),$$
where $r_{j}\in \R$ for $j=1, \dots, M$ and $f_{j, p}(z)$ for $j=1, \dots, M$,
$p=0, \dots, K$ are analytic functions on a disk containing $0$. 
\end{thm}

The convergence and analytic extensions of $q$-traces and pseudo-$q$-traces of 
products of more than two intertwining operators are proved similarly. See 
\cite{H-mod-inv-int}, \cite{F1} and \cite{F2}. 

\renewcommand{\theequation}{\thesection.\arabic{equation}}
\renewcommand{\thethm}{\thesection.\arabic{thm}}
\setcounter{equation}{0} \setcounter{thm}{0}

\section{Convergence results associated to the sewing operations of 
genus-zero Riemann surfaces and determinant lines}

Vertex operator algebras have a geometric definition in terms of 
the partial operad of 
the moduli space of suitable genus-zero Riemann surfaces with punctures and local coordinates and the determinant line bundle over this moduli space. 
See \cite{H-geom-voa}. To prove that a vertex operator algebra 
$V$ indeed satisfies this geometric definition, associated to such 
surfaces, we need to construct 
suitable linear maps from tensor powers of $V$ to the algebraic completion of $V$
using vertex operators (corresponding to 
spheres with three punctures and standard local coordinates) 
and Virasoro operators (corresponding to spheres with two punctures and local coordinates). The main work is to 
prove that these linear maps satisfy some basic properties, including 
a sewing axiom. Since there are conformal anomalies 
(corresponding to central charges for the Virasoro 
operators), the determinant line bundle over the moduli space of such Riemann surfaces
and its powers are also involved. 

The convergence results discussed in the preceding three sections are for
the constructions of correlation functions associated to genus-zero or genus-one
Riemann surfaces with punctures and standard local coordinates 
vanishing at the punctures. To construct and study
correlation functions associated to Riemann surfaces with punctures and general 
local coordinates vanishing at the punctures, we need to 
exponentiate suitable infinite sums of Virasoro operators and prove 
that they behave exactly in the same way as the underlying 
surfaces and determinant lines. The properties that one has to prove include
in particular a convergence involving the exponentials of suitable infinite sums of 
the Virasoro operators and another convergence 
involving the central charge of the Virasoro algebra. These convergence problems related to the Virasoro operators 
were solved by the author in \cite{H-geom-voa}.

For simplicity, we 
use genus-zero Riemann surfaces 
with only two punctures (one positively oriented and 
the other negatively oriented) and local coordinates vanishing at the punctures
to describe these convergence results. 
Such a Riemann surface with punctures and local coordinates is  
conformally equivalent to $\C\cup \{\infty\}$ with the positively oriented puncture
$0$ and negatively oriented puncture $\infty$. Then the local coordinate
vanishing at $0$ becomes a univalent analytic function $f_{0}(w)$
defined near $0$ and vanishing at $0$. Similarly the local coordinate
vanishing at $\infty$ becomes a univalent  analytic  function $f_{\infty}(w)$
defined near $\infty$ and vanishing at $\infty$. 
It is further conformally equivalent to $\C\cup \{\infty\}$ with punctures
$0$ and $\infty$ and with local coordinates given by $f_{0}(w)$ 
and $f_{\infty}(w)$ as above such that the Laurent expansion of $f_{\infty}(w)$
is of the form $\frac{1}{w}+\cdots$, where $\cdots$ are higher order terms
in $\frac{1}{w}$. We denote $\C\cup \{\infty\}$ with such 
punctures and local coordinates by $\Sigma$. 
Such a genus-zero Riemann surface
with two punctures and local coordinates vanishing at the punctures
are said to be canonical. 

Let $[\Sigma]$ be the 
conformal equivalence class of a canonical 
genus-zero Riemann surface $\Sigma$
with two punctures and local coordinates vanishing at the punctures.
As in \cite{H-geom-voa}, we have 
\begin{align*}
f_{0}(w)&=\exp\left(\sum_{j\in Z_{+}}A^{(0)}_{j}w^{j+1}\frac{d}{dw}\right)
a_{0}^{w\frac{d}{dw}}w,\\
f_{\infty}(w)&=\exp\left(\sum_{j\in Z_{+}}A^{(\infty)}_{j}\left(\frac{1}{w}
\right)^{j+1}\frac{d}{d\left(\frac{1}{w}\right)}\right)
\frac{1}{w}
\end{align*}
for some $A^{(0)}_{j}, A^{(\infty)}_{j}\in \C$ for 
$j\in \Z_{+}$ and $a_{0}\in \C^{\times}$. 

Given a vertex operator algebra $V$
(or in general a lower-bounded $\Z$-graded module for the Virasoro algebra),
we define a linear map $\nu_{[\Sigma]}: V\to \overline{V}
=\prod_{n\in \Z}V_{(n)}$ (the algebraic completion 
of $V$) associated to $[\Sigma]$ by 
$$\nu_{[\Sigma]}(v)=\exp\left(-\sum_{j\in Z_{+}}
A^{(\infty)}_{j}L_{V}(-j)\right)
\exp\left(-\sum_{j\in Z_{+}}A^{(0)}_{j}L_{V}(j)\right)
a_{0}^{-L_{V}(0)}v.$$

Let $\Sigma_{1}$ and $\Sigma_{2}$ be two  canonical genus-zero 
Riemann surfaces with two punctures and local coordinates as above. Then 
$\Sigma_{1}$ and $\Sigma_{2}$ are given by the analytic functions
\begin{align*}
f_{0}^{(1)}(w)&
=\exp\left(\sum_{j\in Z_{+}}A_{j}^{(0)}w^{j+1}\frac{d}{dw}\right)
a_{0}^{w\frac{d}{dw}}w,\\
f_{\infty}^{(1)}&
=\exp\left(\sum_{j\in Z_{+}}A^{(\infty)}_{j}\left(\frac{1}{w}
\right)^{j+1}\frac{d}{d\left(\frac{1}{w}\right)}\right)
\frac{1}{w}
\end{align*}
and 
\begin{align*}
f_{0}^{(2)}(w)&
=\exp\left(\sum_{j\in Z_{+}}B_{j}^{(0)}w^{j+1}\frac{d}{dw}\right)
b_{0}^{w\frac{d}{dw}}w,\\
f_{\infty}^{(2)}&
=\exp\left(\sum_{j\in Z_{+}}B^{(\infty)}_{j}\left(\frac{1}{w}
\right)^{j+1}\frac{d}{d\left(\frac{1}{w}\right)}\right)
\frac{1}{w}
\end{align*}
respectively.
Then we have $\nu_{[\Sigma_{1}]}, \nu_{[\Sigma_{2}]}: 
V\to \overline{V}$. We define a series 
$\nu_{[\Sigma_{1}]}\circ \nu_{[\Sigma_{2}]}$
of linear maps from $V$ to $\overline{V}$
by 
$$(\nu_{[\Sigma_{1}]}\circ \nu_{[\Sigma_{2}]})(v)
=\sum_{n\in \Z}\nu_{[\Sigma_{1}]}(\overline{\pi}_{n}
\nu_{[\Sigma_{2}]}(v)),$$
where $\overline{\pi}_{n}$ for $n\in \Z$ is the projection from $\overline{V}$
to $V_{(n)}$. 

On the other hand, if there exists $r\in \R_{+}$ such that 
we can cut disks of radius $r$
from $\Sigma_{1}$ and $\Sigma_{2}$
using the local coordinates vanishing at $0$ on $\Sigma_{1}$ and 
at $\infty$ on $\Sigma_{2}$, respectively, with the remaining 
parts of the surfaces still containing the other punctures,  
we say that $\Sigma_{1}$ can be sewn with $\Sigma_{2}$. 
In this case, 
we identify the boundary of the remaining part of $\Sigma_{1}$
with the boundary of the remaining part of $\Sigma_{2}$ using
the composition of the local coordinate map near $0$ in $\Sigma_{1}$,
the map $w\mapsto \frac{1}{w}$ and the inverse of the 
local coordinate map near $\infty$ in $\Sigma_{2}$ to 
obtain a new genus-zero Riemann 
surfaces with two punctures and local coordinates. We denote it by 
$\Sigma_{1}\;_{1}\infty_{0} \;\Sigma_{2}$. 
The sewing axiom in the geometric definition of vertex operator algebra
states that $\nu_{[\Sigma_{1}]}\circ \nu_{[\Sigma_{2}]}$
is absolutely convergent when $\Sigma_{1}$ can be sewn with $\Sigma_{2}$
and its sum is proportional to 
$\nu_{[\Sigma_{1}\;_{1}\infty_{0} \;\Sigma_{2}]}$. 

To prove the sewing axiom in this case, one first has to prove
the convergence of $\nu_{[\Sigma_{1}]}\circ \nu_{[\Sigma_{2}]}$
when $\Sigma_{1}$ and $\Sigma_{2}$ can be sewn together.
The following result is a special case of a more general 
result proved in \cite{H-geom-voa}:

\begin{thm}[\cite{H-geom-voa}]
Assume that  
$\Sigma_{1}$ can be sewn with $\Sigma_{2}$. 
Then $\nu_{[\Sigma_{1}]}\circ \nu_{[\Sigma_{2}]}$
is absolutely convergent 
in the sense that for $v\in V$ and $v'\in V'$
$$\langle v', (\nu_{[\Sigma_{1}]}\circ \nu_{[\Sigma_{2}]})(v)\rangle
=\sum_{n\in \Z}\langle v', \nu_{[\Sigma_{1}]}(\overline{\pi}_{n}
\nu_{[\Sigma_{2}]}(v))\rangle$$
is absolutely convergent. 
\end{thm}

We briefly explain the idea of the proof of this result.  In fact, 
\begin{align}\label{nu-sewing-1}
&\sum_{n\in \Z}\langle v', \nu_{[\Sigma_{1}]}(\overline{\pi}_{n}
\nu_{[\Sigma_{2}]}(v))\rangle\nn
&\quad=\Biggl\langle v', 
\exp\left(-\sum_{j\in Z_{+}}A^{(\infty)}_{j}L_{V}(-j)\right)
\exp\left(-\sum_{j\in Z_{+}}A^{(0)}_{j}L_{V}(j)\right)
a_{0}^{-L_{V}(0)}\cdot\nn
&\quad\quad \quad\quad 
\cdot \exp\left(-\sum_{j\in Z_{+}}B^{(\infty)}_{j}L_{V}(-j)\right)
\exp\left(-\sum_{j\in Z_{+}}B^{(0)}_{j}L_{V}(j)\right)
b_{0}^{-L_{V}(0)}v\Biggr\rangle,
\end{align}
where the right-hand side should be viewed as a Laurent series in $a_{0}$.
It is proved in \cite{H-geom-voa}
by using formal calculus and properties of the Virasoro operators on $V$
that the right-hand side of (\ref{nu-sewing-1}) is equal to 
\begin{align}\label{nu-sewing-2}
&\sum_{n\in \Z}\Biggl\langle v', 
\exp\left(-\sum_{j\in Z_{+}}A^{(\infty)}_{j}L_{V}(-j)\right)
\exp\left(\sum_{j\in Z_{+}}\Psi_{-j}L_{V}(-j)\right)
\exp\left(\sum_{j\in Z_{+}}\Psi_{j}L_{V}(j)\right)
\cdot\nn
&\quad\quad \quad\quad 
\cdot a_{0}^{-L_{V}(0)}e^{\Psi_{0}L_{V}(0)}e^{\Gamma c}
\exp\left(-\sum_{j\in Z_{+}}B^{(0)}_{j}L_{V}(j)\right)
b_{0}^{-L_{V}(0)}v)\Biggr\rangle,
\end{align}
where $\Psi_{j}$ for $j\in \Z$ and $\Gamma$ are Laurent series 
in $a^{(1)}_{0}$ with polynomials in $A^{(0)}_{j}$ and $B^{(\infty)}_{j}$
for $j\in \Z_{+}$ as coefficients. 
It is proved in \cite{H-geom-voa} 
by using the uniformization theorem and
an old result of Grauert that the Laurent series 
$\Psi_{j}$ for $j\in \Z$ are expansions of analytic functions and therefore
are absolutely convergent. 
It is also proved in \cite{H-geom-voa}
by using the analyticity of 
the canonical isomorphisms between the tensor product of 
the determinant lines of $\Sigma_{1}$ and $\Sigma_{2}$
and the determinant line of $\Sigma_{1}\;_{1}\infty_{0} \;\Sigma_{2}$, 
that the Laurent series $\Gamma$ is the expansion of 
an analytic function and therefore is also absolutely convergent. 
Then (\ref{nu-sewing-2}) and thus the left-hand side of (\ref{nu-sewing-1})
is absolutely convergent.

The convergence result discussed above also applies to lower-bounded 
generalized $V$-modules and even to any lower-bounded graded modules for 
the Virasoro algebra.

The convergence of $\Psi_{j}$ for $j\in \Z$ above was generalized by 
Barron \cite{B1} \cite{B2} to the case of $N=1$ superconformal algebras. 

\renewcommand{\theequation}{\thesection.\arabic{equation}}
\renewcommand{\thethm}{\thesection.\arabic{thm}}
\setcounter{equation}{0} \setcounter{thm}{0}

\section{A higher-genus convergence result of Gui}

Conformal field theories have a geometric formulation given by Kontsevich and 
Segal \cite{G}. Segal in \cite{G} further gave a geometric formulation of chiral conformal field theories (called weakly conformal field theories).
One of the main goal of the mathematical study of conformal field theories
is to construct chiral and full conformal field theories satisfying Segal's axioms.
In particular, one needs to construct correlation functions 
associated to Riemann surfaces with ordered, parametrized and labeled 
boundaries or,
equivalently, Riemann surfaces with
punctures and local coordinates vanishing at punctures, 
from genus-zero Riemann surfaces with one, two 
or three punctures and local coordinates.  
Correlation functions corresponding to the
genus-zero Riemann surface with one puncture and the 
standard local coordinate
are determined by the vacuum of the conformal field theory. 
Correlation functions corresponding to genus-zero Riemann surfaces
with two punctures and 
local coordinates have been discussed in the preceding section and 
are given by the Virasoro operators on modules for 
a vertex operator algebra. Correlation functions corresponding to 
genus-zero Riemann surfaces 
with three 
punctures and standard local coordinates are given by
intertwining operators. 
So to construct chiral conformal field theories, 
one needs to construct correlation functions associated to 
arbitrary Riemann surfaces with punctures and local coordinates 
from the vacuum, the Virasoro operators and intertwining operators. 
The genus-zero correlation functions and genus-one correlation functions 
were constructed in \cite{H-diff-eqn} and \cite{H-mod-inv-int}, 
respectively (see also Sections 3 and 4 for the convergence problems
associated to these constructions). 

To construct higher-genus chiral correlation functions from 
genus-zero and genus-one chiral correlation functions, we need to 
prove a higher-genus convergence. For rational conformal 
field theories, this was in fact 
stated as a conjecture 
in \cite{Z1} and \cite{H-open-prob}.
This conjecture was proved in 2020 by Gui \cite{G}. 

We now describe this higher-genus convergence result. 
Let $W_{1}, \dots, W_{n}$ be grading-restricted generalized $V$-modules.
For $w_{1}\in W_{1}$, $\dots$, $w_{n}\in W_{n}$,
an $n$-point genus-$g$ correlation function associated to $W_{1},
\dots, W_{n}$  is a linear map from 
$W_{1}\otimes \cdots \otimes W_{n}$ to the space of 
multivalued analytic functions 
on the moduli space of genus-$g$
Riemann surfaces with $n$ punctures and local coordinates
vanishing at the punctures satisfying suitable conditions.
Here for simplicity we omit the description of these conditions. 

Let $\Sigma_{1}$  be a genus-$g_{1}$
Riemann surface with 
$n_{1}$ punctures and local coordinates vanishing at the punctures
and $\Sigma_{2}$ a genus-$g_{2}$
Riemann surface with 
$n_{2}$ punctures and local coordinates vanishing at the punctures.
If there exists $r\in \R_{+}$ such that 
we can cut disks of radius $r$
from $\Sigma_{1}$ and $\Sigma_{2}$
using the local coordinates vanishing at the $i$-th puncture on $\Sigma_{1}$ 
and at the $j$-th puncture on $\Sigma_{2}$, respectively, with the remaining 
parts of the surfaces still containing the other punctures,  
we say that $\Sigma_{1}$ can be sewn with $\Sigma_{2}$ at the $i$-th 
puncture on $\Sigma_{1}$ and the $j$-th puncture on $\Sigma_{2}$. 
In this case, 
we can identify the boundary of the remaining part of $\Sigma_{1}$
with the boundary of the remaining part of $\Sigma_{2}$ using
the composition of the local coordinate map near the $i$-th puncture 
on $\Sigma_{1}$,
the map $w\mapsto \frac{1}{w}$ and the inverse of the 
local coordinate map near the $j$-th puncture on $\Sigma_{2}$ to 
obatin a new Riemann 
surfaces with punctures and local coordinates. 

Let $\psi_{1}$ be an $n_{1}$-point genus-$g_{1}$ correlation function 
associated to $W_{1}, \dots, W_{n_{1}}$
and $\psi_{2}$ an $n_{2}$-point genus-$g_{2}$-correlation function
associated to $\widetilde{W}_{1}, \dots, \widetilde{W}_{n_{2}}$.
Assume that $\widetilde{W}_{j}=W_{i}'$. 
One axiom for chiral conformal field theories requires
that the series 
\begin{align}\label{sewing-of-2}
&\sum_{k\in \Z_{+}}(\psi_{1}(w_{1}\otimes \cdots w_{i-1}
\otimes w_{i}^{(k)}\otimes w_{i+1}
\otimes \cdots \otimes w_{n_{1}}))([\Sigma_{1}])\cdot\nn
&\quad\quad\quad\quad\cdot 
(\psi_{2}(\tilde{w}_{1}\otimes \cdots \tilde{w}_{j-1}
\otimes (w_{i}^{(k)})'\otimes \tilde{w}_{j+1}
\otimes \cdots \otimes \tilde{w}_{n_{2}}))([\Sigma_{2}])
\end{align}
be absolutely convergent when 
$\Sigma_{1}$ can be sewn with $\Sigma_{2}$ at the $i$-th 
puncture on $\Sigma_{1}$ and the $j$-th puncture on $\Sigma_{2}$,
where $\{w_{i}^{(k)}\}_{k\in \Z_{+}}$ 
and $\{(w_{i}^{(k)})'\}_{k\in \Z_{+}}$ are dual homogeneous basis of 
$W_{i}$ and $W_{i}'=\widetilde{W}_{j}$. 
This is the higher-genus convergence problem for the sewing 
of two Riemann surfaces. The convergence of 
products of intertwining operators discussed in Section 3 is the special case
of this convergence. 

There is also another convergence problem for the self sewing 
of one Riemann surface. 
Let $\Sigma$ be a genus-$g$ Riemann surface with $n$ punctures and local coordinates vanishing at punctures. 
If there exists $r\in \R_{+}$ such that 
we can cut disks of radius $r$
from $\Sigma$ 
using the local coordinates vanishing at the $i$-th and the $j$-th 
punctures on $\Sigma$ with the remaining 
parts of the surfaces still containing the other punctures,  
we say that $\Sigma$ can be sewn at the $i$-th 
 and $j$-th punctures. In this case, we can also obtain a new 
 genus-$g+1$ Riemann surface with $n-2$ punctures and 
 local coordinates vanishing at the punctures by sewing 
 $\Sigma$ at the $i$-th and $j$-th punctures using the same 
 procedure as in the case of two surfaces above. 
 
Let $\psi$ be an $n$-point genus-$g$ correlation function 
associated to $W_{1}, \dots, W_{n}$. Assume that 
$W_{j}=W_{i}'$. Assume that $i<j$. 
Then one axiom for chiral conformal field theories requires
that the series 
\begin{align}\label{self-sewing}
&\sum_{k\in \Z_{+}}(\psi(w_{1}\otimes \cdots w_{i-1}
\otimes w_{i}^{(k)}\otimes w_{i+1}\otimes \cdots
\otimes w_{j-1}
\otimes (w_{i}^{(k)})'\otimes w_{j+1}
\otimes \cdots \otimes w_{n}))([\Sigma])
\end{align}
be absolutely convergent when 
$\Sigma$ can be sewn at the $i$-th 
 and $j$-th punctures. 
This is the convergence problem for the self sewing 
of one Riemann surface. 
The convergence of the $q$-traces
of products of geometrically-modified intertwining operators
discussed in Section 4
is the special case of this convergence. 

The following theorem is proved by Gui in \cite{G}:

\begin{thm}[Gui \cite{G}]\label{higher-genus}
Let $V$ be a $C_{2}$-cofinite 
vertex operator algebra containing no nonzero 
elements of negative weights. If the grading-restricted generalized
$V$-modules involved are 
finitely generated, then 
(\ref{sewing-of-2}) and (\ref{self-sewing}) are absolutely convergent
when 
$\Sigma_{1}$ can be sewn with $\Sigma_{2}$ at the $i$-th 
puncture on $\Sigma_{1}$ and the $j$-th puncture on $\Sigma_{2}$
and when 
$\Sigma$ can be sewn at the $i$-th 
 and $j$-th punctures, respectively. 
 \end{thm}
 
As in the proof of the convergence of products 
of intertwining operators and the proof of the 
convergence of $q$-traces and pseudo-$q$-traces of products of
geometrically-modified
intertwining operators
discussed Sections 3 and 4, respectively,
 this theorem is proved in \cite{G} by deriving differential equations. 
 But in this case the derivation of the differential equations involving 
 analytic functions on the moduli space of higher-genus Riemann surfaces.
These analytic functions are much more difficult to study than those
on the moduli spaces of genus-zero and genus-one Riemann surfaces. 
For example, the differential equations satisfied by 
$q$-traces and pseudo-$q$-traces of products of
intertwining operators were derived using the $q$-expansions of the 
derivatives of the Weierstrass function (see \cite{H-mod-inv-int}).
We need similar results for functions on the moduli space 
of higher-genus surfaces. 
This difficult in the higher-genus case was overcame in \cite{G} 
by using a theorem of Grauert in complex analysis. 

\renewcommand{\theequation}{\thesection.\arabic{equation}}
\renewcommand{\thethm}{\thesection.\arabic{thm}}
\setcounter{equation}{0} \setcounter{thm}{0}

\section{Convergence conjectures and problems
 in orbifold conformal field theory}	

Orbifold conformal field theories are 
conformal field theories constructed from 
known conformal field theories and their automorphisms. 
In the framework of the representation theory of vertex operator algebras, 
orbifold conformal field theory is the 
study of twisted intertwining operators among 
(generalized) twisted modules. In this section, we discuss the 
convergence conjectures and problems for orbifold conformal field theories. 

Let $V$ be a vertex operator algebra and $g$ an automorphism of $V$. 
A lower-bounded generalized $g$-twisted module is a $\C$-graded vector space 
$W=\coprod_{n\in \C}W_{[n]}$ 
such that $W_{[n]}=0$ when $\Re(n)$ is sufficiently negative,
equipped with a
twisted vertex operator map 
\begin{align*}
Y_{W}^{g}: V\otimes W&\to W\{z\}[\log z]\nn
v\otimes w&\mapsto Y_{W}^{g}(v, z)w
\end{align*}
satisfying suitable axioms, including in particular an 
equivariance property and a duality property which 
requires that products of twisted vertex operators are convergent 
in suitable regions and the associativity and commutativity for 
twisted intertwining operators hold.
To construct lower-bounded generalized twisted modules,  
the convergence of of products of twisted vertex operators
can be proved using the method in Section 2 (see
 \cite{H-const-twisted-mod}). 
In \cite{DLM}, by using the method of Zhu \cite{Z},
Dong, Li and Mason generalized the convergence and analytic extension 
results of Zhu \cite{Z} to the convergence and analytic extension
of $q$-traces of twisted vertex operators on a 
$g$-twisted module associated to a finite order automorphism $g$ of 
a $C_{2}$-cofinite vertex operator algebra $V$.
But orbifold conformal field theories
are about twisted intertwining operators among twisted modules.
In general, the convergences of products and (pseudo-)q-traces of 
products of twisted intertwining operators are still conjectures. 

Let $W_{1}$, $W_{2}$ and $W_{3}$ be generalized $g_{1}$-,
$g_{2}$- and $g_{3}$-twisted $V$-modules, respectively.
A twisted intertwining operator $\Y$ of type $\binom{W_{3}}{W_{1}W_{2}}$ 
is a linear map 
\begin{align*}
\Y: W_{1}\otimes W_{2}&\to W_{3}\{z\}[\log z]\nn
w_{1}\otimes w_{2}&\mapsto \Y(w_{1}, z)w_{2}
\end{align*}
satisfying a duality property and an $L(-1)$-derivative property. 
In particular, just as intertwining operators among (untwisted) 
generalized $V$-modules, for $w_{1}\in W_{1}$, we have 
 $$\Y(w_{1}, z)\in \hom(W_{3}, W_{2})[\log z]\{z\}.$$ 

As in the case of intertwining operators among (untwisted) 
generalized $V$-modules, let 
$W_{1}$, $W_{2}$, $W_{3}$, $W_{4}$ and $W_{5}$ 
 be generalized $g_{1}$-, $g_{2}$-, $g_{3}$-, $g_{3}$- and 
 $g_{5}$-twisted $V$-modules, respectively, 
  and $\Y_{1}$ and $\Y_{2}$ twisted
 intertwining operators of types $\binom{W_{4}}{W_{1}W_{5}}$ 
 and $\binom{W_{5}}{W_{2}W_{3}}$,
 respectively. Then for $w_{1}\in W_{1}$ and $w_{2}\in W_{2}$, 
 $$\Y_{1}(w_{1}, z_{1})
 \Y_{1}(w_{2}, z_{2})\in \hom(W_{3}, W_{4})[\log z_{2}, \log z_{2}]
 \{z_{1}, z_{2}\}.$$
We have a series in $\C$
 \begin{equation}\label{twisted-int-op-prod}
 \langle w_{4}', \Y_{1}(w_{1}, z_{1})
 \Y_{1}(w_{2}, z_{2})w_{3}\rangle\lbar_{\log z_{1}=l_{p}(z_{1}),\; \log z_{2}=l_{p}(z_{2})}
 \end{equation}
 for $w_{1}\in W_{1}$ and $w_{2}\in W_{2}$, 
 $w_{3}\in W_{3}$ and $w_{4}'\in W_{4}'$.

Note that $W_{1}$, $W_{2}$, $W_{3}$, $W_{4}$ and $W_{5}$
are generalized $V^{G}$-modules where $G$ is the fixed point subalgebra
of $V$ under the group $G$ generated by 
$g_{1}, g_{2}, g_{3}, g_{4}, g_{5}$. Also $\Y_{1}$ and $\Y_{2}$
are intertwining operators 
of types $\binom{W_{4}}{W_{1}W_{5}}$ 
 and $\binom{W_{5}}{W_{2}W_{3}}$ when 
 $W_{1}$, $W_{2}$, $W_{3}$, $W_{4}$ and $W_{5}$ are viewed as 
 generalized $V^{G}$-modules. Thus if 
 $W_{1}$, $W_{2}$, $W_{3}$ and  $W_{4}'$
 are quasi-finite-dimensional and $C_{1}$-cofinite as $V^{G}$-modules,
by Theorem \ref{conv-prod-int}, 
 (\ref{twisted-int-op-prod}) 
 is absolutely convergent in the region $|z_{1}|>|z_{2}|>0$ and 
 can be analytic extended as in Theorem 
\ref{conv-prod-int}. This approach indeed works
in the case that $V$ satisfies the three conditions in Theorem
\ref{higher-genus} and $G$ is a finite solvable group because
Canahan and Miyamoto proved in \cite{CM} that in this case, $V^{G}$
also satisfies these three conditions (see also \cite{Mc} for a new proof
of this result). 

But for an infinite group or a finite nonsolvable
group $G$, $V^{G}$ being $C_{2}$-cofinite 
or even $C_{1}$-cofinite when $V$ is 
$C_{2}$-cofinite is still an open problem. 
Instead of trying to prove $V^{G}$ satisfies the conditions 
needed, the author proposed in \cite{H-open-prob}
and \cite{H-orbifold-theory}
a program to study 
orbifold conformal field theories by studying directly 
twisted intertwining operators 
among suitable generalized twisted $V$-modules. 
In this program, we need to prove in particular that 
(\ref{twisted-int-op-prod}) 
 is absolutely convergent in the region $|z_{1}|>|z_{2}|>0$
and the sum can be analytic extended as in Theorem 
\ref{conv-prod-int}. 
As in the case of  intertwining operators among (untwisted) 
generalized $V$-modules, we expect that 
(\ref{twisted-int-op-prod}) is absolutely convergent only when 
$W_{1}$, $W_{2}$, $W_{3}$, $W_{4}$ and $W_{5}$
as generalized twisted $V$-modules (not as generalized $V^{G}$-modules)
satisfy certain conditions. 

Though (\ref{twisted-int-op-prod}) looks completely the same as 
 (\ref{int-op-prod}), it is much more difficult to study
 since the twisted vertex operators 
 for generalized twisted modules in general involve nonintegral
 powers and logarithms of the variables. It is especially difficult 
 to study in the case that the automorphisms 
 $g_{1}, g_{2}, g_{3}, g_{4}, g_{5}$ do not commute with each other. 
There is still no general convergence result yet. 
But we have the following conjecture formulated in \cite{H-open-prob}
and \cite{H-orbifold-theory}:

\begin{conj}\label{conj-orbi-prod--tw-int}
Let $V$ be a vertex operator algebra satisfying the following conditions:
(i) $V$ is of positive energy, that is, $V_{(n)}=0$ for $n<0$ 
and $V_{(0)}=\C\one$, and $V$ is equivalent to $V'$ as a $V$-module. 
(ii) $V$ is $C_{2}$-cofinite. (iii) Every grading-restricted generalized 
$V$-module is completely reducible. Let $G$ be a finite group 
of automorphisms of $V$. Let  
$g_{1}, g_{2}, g_{3}, g_{4}, g_{5}\in G$, 
$W_{1}$, $W_{2}$, $W_{3}$, $W_{4}$ and $W_{5}$ 
grading-restricted generalized $g_{1}$-, $g_{2}$-, $g_{3}$-, $g_{4}$- and 
 $g_{5}$-twisted $V$-modules, respectively, 
  and $\Y_{1}$ and $\Y_{2}$ twisted
 intertwining operators of types $\binom{W_{4}}{W_{1}W_{5}}$ 
 and $\binom{W_{5}}{W_{2}W_{3}}$,
 respectively. Then for $w_{1}\in W_{1}, w_{2}\in W_{2}$, 
 $w_{3}\in W_{3}$ and $w_{4}'\in W_{4}'$, 
the series (\ref{twisted-int-op-prod}) is absolutely convergent in the 
region $|z_{1}|>|z_{2}|>0$ and its sum can be analytically 
continued to a multivalued analytic function on $M^{2}$. 
\end{conj}

The general form of this convergence conjecture is for the product
of more than two intertwining operators. See \cite{H-orbifold-theory}
for details. 

For orbifold conformal field theories, we also need to 
study $q$-traces or pseudo-$q$-traces of products of intertwining operators.
Let $g_{1}, g_{2}, g_{3}, g_{4}$ be automorphism of $V$, 
$W_{1}$, $W_{2}$, $W_{3}$ and $W_{4}$ 
grading-restricted generalized $g_{1}$-, $g_{2}$-, $g_{3}$- and 
 $g_{4}$-twisted $V$-modules, respectively, 
  and $\Y_{1}$ and $\Y_{2}$ twisted
 intertwining operators of types $\binom{W_{3}}{W_{1}W_{4}}$ 
 and $\binom{W_{4}}{W_{2}W_{3}}$,
 respectively. Let $P$ be a finite-dimensional associative algebra
 and $\phi$ a symmetric linear function on $P$. Assume that 
$W_{3}$ is a projective right $P$-module such that 
its twisted vertex operators  and 
$$\mathcal{Y}_{1}(\mathcal{U}_{W_{1}}(q_{z_{1}})w_{1}, 
q_{z_{1}})
\mathcal{Y}_{2}(\mathcal{U}_{W_{2}}(q_{z_{2}})w_{2}, q_{z_{2}})$$
commute with the action of $P$ on $W_{3}$. Then we have 
the pseudo-$q$-trace
\begin{align}\label{pseudo-q-trace-tw}
&\tr^{\phi}_{W_{3}}\mathcal{Y}_{1}(\mathcal{U}_{W_{1}}(q_{z_{1}})w_{1}, 
q_{z_{1}})
\mathcal{Y}_{2}(\mathcal{U}_{W_{2}}(q_{z_{2}})w_{2}, q_{z_{2}})
q^{L(0)-\frac{c}{24}}\nn
&\quad =\sum_{n\in \C}\left(\tr^{\phi}_{(W_{3})_{[n]}}
\pi_{n}\mathcal{Y}_{1}(\mathcal{U}_{W_{1}}(q_{z_{1}})w_{1}, 
q_{z_{1}})
\mathcal{Y}_{2}(\mathcal{U}_{W_{2}}(q_{z_{2}})w_{2}, q_{z_{2}})
q^{L_{W_{3}}(0)-\frac{c}{24}}\lbar_{(W_{3})_{[n]}}
\right)\nn
&\quad =\sum_{n\in \C}\left(\tr^{\phi}_{(W_{3})_{[n]}}
\pi_{n}\mathcal{Y}_{1}(\mathcal{U}_{W_{1}}(q_{z_{1}})w_{1}, 
q_{z_{1}})
\mathcal{Y}_{2}(\mathcal{U}_{W_{2}}(q_{z_{2}})w_{2}, q_{z_{2}})
q^{n-\frac{c}{24}}e^{(\log q) L_{W}(0)_{N}}\lbar_{(W_{3})_{[n]}}\right).
\end{align}

We have the following convergence conjecture on 
pseudo-$q$-traces of products of intertwining operators
given in \cite{H-open-prob}
and \cite{H-orbifold-theory}:

\begin{conj}
Let $V$ be a vertex operator algebra satisfying the conditions in 
Conjecture \ref{conj-orbi-prod--tw-int} and $G$ a finite group 
of automorphisms of $V$. Then for $g_{1}, g_{2}, g_{3}, g_{4}\in G$, 
the series (\ref{pseudo-q-trace-tw}) with $q=q_{\tau}=e^{2\pi i\tau}$
is absolutely convergent in the region $1>|q_{z_{1}}|
>|q_{z_{2}}|>|q_{\tau}|>0$ and can be analytically extended
to a multivalued analytic function in the region given by
$\Im(\tau)>0$ (here $\Im(\tau)$ is the imaginary part of $\tau$), 
$z_{1}\ne z_{2}+k\tau+l$ for $k, l\in \mathbb{Z}$. Moreover, the singular point
$z_{1}=z_{2}+k\tau+l$ for each $k, l\in \mathbb{Z}$ is regular, that is,
any branch of the multivalued analytic function can be expanded 
in a neighborhood of 
the singular point $z_{1}=z_{2}+k\tau+l$ as a series of the form
$$\sum_{p=0}^{K}
\sum_{j=1}^{M}(z_{1}-z_{2}+k\tau+l)^{r_{j}}(\log (z_{1}-z_{2}+k\tau+l))^{p}
f_{j, p}(z_{1}-z_{2}+k\tau+l),$$
where $r_{j}\in \R$ for $j=1, \dots, M$ and $f_{j, p}(z)$ for $j=1, \dots, M$,
$p=0, \dots, K$ are analytic functions on a disk containing $0$. 
\end{conj}

The general form of this convergence conjecture is for pseudo-$q$-traces 
of products
of more than two intertwining operators. See \cite{H-orbifold-theory}
for details. 

In the general case, we have the following problem:

\begin{prob}
Let $V$ be  a vertex operator algebra and let $G$ be a group of automorphisms of $V$. Under what conditions do products and pseudo-$q$-traces 
of products of 
twisted intertwining operators 
among the grading-restricted generalized $g$-twisted $V$-modules 
for $g\in G$ converge and have analytic extensions as in the 
conjectures above? 
\end{prob}

\renewcommand{\theequation}{\thesection.\arabic{equation}}
\renewcommand{\thethm}{\thesection.\arabic{thm}}
\setcounter{equation}{0} \setcounter{thm}{0}

\section{Convergence problems in the cohomology theory of vertex algebras}

In \cite{H-cohomology} and \cite{H-1st-2nd-coh},
the author introduced a cohomology theory for grading-restricted vertex 
algebras 
and showed that the cohomology for a grading-restricted 
vertex algebra has the properties that a conhomology theory must have.
Let $V$ be a grading-restricted vertex algebra and $W$ 
a grading-restricted generalized $V$-module. Let 
$\overline{W}=\prod_{n\in \C}W_{[n]}$ be the algebraic completion of 
$W$. A $\overline{W}$-valued rational function in $z_{1}, \dots, z_{n}$
is a $\overline{W}$-valued function  $f(z_{1}, \dots, z_{n})$ such that 
$\langle w', f(z_{1}, \dots, z_{n})\rangle$ is a rational function in $z_{1}, 
\dots, z_{n}$ for $w'\in W'$. 
In this cohomology theory, $n$-cochains with coefficients in $W$
are maps from the $n$-th tensor power of $V$ to the space of 
$\overline{W}$-valued rational functions in variables 
$z_{1}, \dots, z_{n}$ with the only possible poles 
$z_{i}-z_{j}=0$ for $i\ne j$, satisfying several conditions, including in particular 
a condition that the series obtained by composing these maps 
with vertex operator maps
for $V$ and for $W$  involving additional $m$ variables  
$z_{n+1}, \dots, z_{m+n}$ are absolutely convergent 
in suitable regions and can be analytically extended to rational functions 
in $z_{1}, \dots, z_{m+n}$ with the only possible poles at 
$z_{i}=z_{j}$ for $i\ne j$, $i, j=1, \dots, m+n$. 
This convergence is important since the coboundary operator 
is defined using the rational functions obtained from these convergent series.

Since cochains in this cohomology theory by definition must satisfy 
such a convergence condition, results and explicit calculations 
in this cohomology theory are always based on some basic 
convergence results or assumptions. Though the series involved 
should be convergent to rational functions, the method in Section 2 
cannot be applied because cochains does not satisfy properties such as 
weak commutativity or weak associativity. To understand this cohomology theory
and apply it to solve mathematical 
problems, we need to find algebraic conditions 
on the vertex algebra and modules and to develop new techniques to 
prove this type of convergence under these algebraic conditions. 

It is proved in \cite{H-1st-2nd-coh} that $1$-cochains always satisfy 
the convergence condition. So here we use $2$-cochain to 
discuss the convergence condition. Given a grading-restricted vertex algebra
$V$ and a grading-restricted (or lower-bounded) generalized 
$V$-module $W$, a $2$-cochain composable with 
$m$ vertex operators is equivalent to a linear map 
\begin{align}\label{2-cochain}
\Psi: V\otimes V&\to W((x))\nn
v_{1}\otimes v_{2}&\mapsto \Psi(v_{1}, x)v_{2}
\end{align}
satisfying certain conditions, including in particular, 
the condition that the series obtained from 
the compositions of $\Psi$ with $m$ vertex operators 
are absolutely convergent in suitable regions. For example, 
for $w'\in W'$, $v, v_{1}, \dots, v_{k+l+1}\in V$, 
$$\langle w', Y_{W}(v_{1}, z_{1})\cdots Y_{W}(v_{k}, z_{k})
 \Psi(v, z)
Y_{V}(v_{k+1}, z_{k+1})\cdots Y_{V}(v_{k+l}, z_{k+l})
v_{k+l+1}\rangle$$
should be absolutely convergent in the region 
$|z_{1}|>\cdots |z_{k}|>|z|>|z_{k+1}|>\cdots >|z_{k+l}|>0$
to a rational function in $z_{1}, \dots, z_{k+l}, z$ with the only possible 
poles $z_{i}=0$ for $i=1, \dots, k+l$, 
$z=0$, $z_{i}=z_{j}$ for $1\le i<j\le k+l$ and $z_{i}=z$ for 
$i=1, \dots, k+l$.
Another type of compositions is given by iterates, for example, 
$$\langle w', \Psi(Y_{V}(v_{1}, z_{1}-z_{m+1})\cdots 
Y_{V}(v_{m}, z_{m}-z_{m+1})v_{m+1}, z_{m+1})v_{m+2}\rangle.$$
This series is required to be absolutely convergent in the region 
$|z_{m+1}|>|z_{1}-z_{m+1}|>\cdots >|z_{m}-z_{m+1}|>0$
to a rational function in $z_{1}, \dots, z_{m+1}$ of the form above. 
There are certainly many different ways to compose 
$\Psi$ with $m$ vertex operators. They are all required to be 
absolutely convergent  in suitable regions to 
rational functions in $z_{1}, \dots, z_{m+1}$ of the form above.

To calculate explicitly the cohomology of a grading-restricted vertex algebra, 
we need to find all the cochains first. Thus the first problem in 
such a calculation is to determine all the maps 
from the $n$-th tensor power of $V$ to the space of 
$\overline{W}$-valued rational functions in variables 
$z_{1}, \dots, z_{n}$ satisfying the convergence condition. 
This is in general not an easy problem, even for a relatively 
simple vertex algebra. For example, if we want to calculate 
the second cohomology of a grading-restricted vertex algebra $V$,
we need to determine all those maps of the form (\ref{2-cochain})
such that the convergence condition holds. For example, in the case $m=1$
above, we need to determine in particular whether 
\begin{align*}
&\langle w', \Psi(v_{1}, z_{1})Y_{V}(v_{2}, z_{2})v_{3}\rangle,\\
&\langle w', Y_{W}(v_{1}, z_{1})\Psi(v_{2}, z_{2})v_{3}\rangle
\end{align*}
are absolutely convergent in the region $|z_{1}|>|z_{2}|>0$ 
to a rational function in $z_{1}$ and $z_{2}$ with the only possible poles at 
$z_{1}=0$, $z_{2}=0$ and $z_{1}-z_{2}=0$. 
Note that $\Psi(v_{1}, z_{1})$ and $Y_{V}(v_{2}, z_{2})$
or $Y_{W}(v_{1}, z_{1})$ and $\Psi(v_{2}, z_{2})$ do not have to satisfy 
the weak commutativity, as we have discussed in Section 2, 
the method in Section 2 cannot be used to determine such 
$\Psi$. So we have the following problem:

\begin{prob}
Is there a general method that can be used to determine all the cochains?
Are there algebraic conditions on $V$ and $W$ such that 
we can determine all the cochains coposable with 
$m$ vertex operators using these algebraic conditions. 
\end{prob}

 Another convergence problem related to the cohomology 
 theory appeared in the work \cite{HQ}
 of Qi and the author. It has been proved in \cite{H-1st-2nd-coh}
 that the first cohomology of a grading-restricted vertex algebra
 $V$ with coefficients in a  grading-restricted generalized $V$-module
 $W$ is isomorphic to the space of derivations from $V$ to $W$ 
 modulo the space of inner derivations. 
 In \cite{HQ}, it is proved that if every derivation from $V$ to $W$ is 
 the sum of an inner derivation and a derivation called zero-mode derivation
 for every  $\Z$-graded bimodules when $V$ is viewed as a 
 meromorphic open-string vertex algebra, then every 
 lower-bounded generalized $V$-module satisfying a composability 
 condition is completely reducible. In this result, the complete reducibility 
 holds only for  lower-bounded generalized $V$-module satisfying a 
 composability  condition. This composability condition is in
 fact a convergence condition. 
 
 We now describe this composability condition.
 Let $W$ be a lower-bounded generalized $V$-module and 
$W_{2}$ a $V$-submodule  of $W$. 
We say that the pair $(W, W_{2})$ satisfies the composability condition
if there exists a graded subspace $W_{1}$ of $W$ such that 
$W=W_{1}\oplus W_{2}$ as a graded 
vector space and  
\begin{equation}\label{comp-prod}
\langle w_{2}', Y_{W_{2}}(v_{1}, z_{1})\cdots Y_{W_{2}}(v_{k}, z_{k})\pi_{W_{2}}Y_{W}(v, z)\pi_{W_{1}}
Y_{W_{1}}(v_{k+1}, z_{k+1})\cdots Y_{W_{1}}(v_{k+l}, z_{k+l})
w_{1}\rangle
\end{equation}
for $v, v_{1}, \dots, v_{k+l}\in V$, $w_{2}'\in W_{2}'$ and $w_{1}\in W_{1}$
is absolutely convergent in the region 
$|z_{1}|>\cdots |z_{k}|>|z|>|z_{k+1}|>\cdots >|z_{k+l}|>0$
to a rational function in $z_{1}, \dots, z_{k+l}, z$ with the only possible 
poles $z_{i}=0$ for $i=1, \dots, k+l$, 
$z=0$, $z_{i}=z_{j}$ for $1\le i<j\le k+l$ and $z_{i}=z$ for 
$i=1, \dots, k+l$ such that the orders of the poles satisfy 
some conditions which we omit here. If for every 
proper nonzero left $V$-submodule $W_{2}$ of $W$, 
the pair $(W, W_{2})$ satisfies the composability condition, 
we say that $W$ satisfies the composability condition. 

Assume that $V$ contains a subalgebra $V_{0}$ such that the 
following conditions for intertwining operators among 
grading-restricted generalized $V_{0}$-modules are satisfied:

\begin{enumerate}

\item  For any $n\in \Z_{+}$, 
products of $n$ intertwining operators among grading-restricted 
generalized $V_{0}$ modules 
evaluated at $z_{1}, \dots, z_{n}$
are absolutely convergent in the region $|z_{1}|>\cdots >|z_{n}|>0$ and
can be analytically extended to  (possibly multivalued) analytic functions in $z_{1}, \dots, z_{n}$
with the only possible 
singularities (branch points or poles)
$z_{i}=0$ for $i=1, \dots, n$ and $z_{i}=z_{j}$ for 
$i, j=1, \dots, n$, $i\ne j$. 

\item The associativity of intertwining operators
among grading-restricted generalized $V_{0}$-modules holds. 

\end{enumerate}
Then it is proved in \cite{HQ} that for a lower-bounded generalized 
$V$-module $W$, a lower-bounded generalized $V$-submodule 
$W_{2}$ of $W$ and a lower-bounded generalized $V_{0}$-submodule 
$W_{1}$ of $W$ such that $W=W_{1}\oplus W_{2}$, 
the pair $(W, W_{2})$ satisfies the composability condition. In particular,
$W$ satisfies the composability condition. 

This result on the composability requires that there is a nice 
subalgebra $V_{0}$ of $V$. This is in general not true. 

In \cite{Q}, Qi studied the composability condition in the case that 
$V$ is a Virasoro vertex operator algebra. For certain special $V$-modules
and their submodules, he proved the composability conditions. 

We have the following main convergence problem related to this 
complete reducibility theorem in \cite{HQ}:

\begin{prob}
Are there algebraic conditions on $V$, $W$ and $W_{2}$ such that 
the pair $(W, W_{2})$ satisfies the composability condition if these 
algebraic conditions are satisfied?
\end{prob}

\noindent {\small \sc Department of Mathematics, Rutgers University,
110 Frelinghuysen Rd., Piscataway, NJ 08854-8019, USA}

\noindent {\em E-mail address}: {\tt yzhuang@math.rutgers.edu}

\end{document}